\documentclass[12pt]{article}
\usepackage{amsmath,amssymb}
\usepackage{theorem}
\usepackage{algorithm}
\usepackage{graphicx}
\usepackage[usenames]{color}
\usepackage{multirow}

\newcommand{\mybf}[1]{\boldsymbol{#1}}

  {\end{list}}

\makeatletter

  \gdef\listctr{list\romannumeral\the\@listdepth}\expandafter
  
\makeatother

  {\end{list}}

\title{On the filtering effect of iterative regularization algorithms for linear least-squares problems}

\author{A Cornelio, F Porta, M Prato and L Zanni\\
\small Dipartimento di Scienze Fisiche, Informatiche e Matematiche\\
\small Universit\`a di Modena e Reggio Emilia\\ 
\small Via Campi 213/b, 41125 Modena, Italy\\
\small \texttt{anastasia.cornelio@unimore.it, federica.porta@unimore.it}\\
\small \texttt{marco.prato@unimore.it, luca.zanni@unimore.it}}

\begin{document}
\maketitle

\begin{abstract}
Many real-world applications are addressed through a linear least-squares problem formulation, whose solution is calculated 
by means of an iterative approach. A huge amount of studies has been carried out in the optimization field to provide the fastest 
methods for the reconstruction of the solution, involving choices of adaptive parameters and scaling matrices. However, in presence 
of an ill-conditioned model and real data, the need of a regularized solution instead of the least-squares one changed the point 
of view in favour of iterative algorithms able to combine a fast execution with a stable behaviour with respect to the restoration 
error. In this paper we want to analyze some classical and recent gradient approaches for the linear least-squares problem by 
looking at their way of filtering the singular values, showing in particular the effects of scaling matrices and non-negative 
constraints in recovering the correct filters of the solution.
\end{abstract}

\section{Introduction}

We consider the classical linear least-squares problem
\begin{equation}\label{eq1}
\min_{\mybf{x} \in \mathbb{R}^n} f(\mybf{x}) \equiv \frac{1}{2} \| \mybf{A} \mybf{x} - \mybf{b} \|_2^2,
\end{equation}
where $\mybf{A}$ is a full-rank $m \times n$ matrix, with $m \geq n$, and $\mybf{b} \in \mathbb{R}^m$. 
Such issue is typically addressed in presence of an ill-posed inverse problem 
\begin{equation}\label{eq2}
\mybf{A} \mybf{x} = \mybf{b}, 
\end{equation}
since the existence and uniqueness of the least-squares solution $\mybf{x}^\dagger$ together with its continuously data dependence 
(that is guaranteed in the discrete case) lead to a well-posed problem. However, especially in the case in which the linear system 
\eqref{eq2} arises from the discretization of a continuous ill-posed inverse problem, the switch to the least-squares problem 
\eqref{eq1} does not avoid the ill-conditioning pathology, that amplifies the noise affecting the data with the result of 
a meaningless reconstructed solution. A classical example occurs in the image deblurring problem \cite{BER98,HAN06}, in which 
$\mybf{b}$ is a blurred and noisy version of an unknown image $\mybf{x}$ and $\mybf{A}$ describes the transformation from the 
target to the measured data values. The numerical instability provided by the ill-conditioning can be countered by looking for 
a regularized solution of \eqref{eq1}, obtained by means of an approximation of the original problem with a one-parameter family 
of better conditioned ones. A strategy for the choice of the ``best'' parameter is finally needed, according to some criteria 
based for example on the knowledge of the noise level or the analysis of the residuals \cite{ENG00,HAN97,VOG02}. Besides the 
classical direct regularization approaches as the truncated singular value decomposition or the Tikhonov method, an increasing 
interest has been devoted to iterative regularization strategies, that are in general more computationally effective for large 
scale problems and allow an easier introduction of constraints on the desired solution (e.g., non-negativity, flux conservation) 
\cite{BAR06,BEN10,BER98,BER08,VOG02}. Although in general the regularization property of iterative approaches is not 
theoretically proved, most of them exhibits the so-called {\em semiconvergence effect}, i.e., the reconstruction error decreases 
until one optimal iterate and then diverges. Such behaviour has been analyzed e.g. by Nagy \& Palmer \cite{NAG03} and 
Donatelli \& Serra Capizzano \cite{DON08} in terms of {\em filter factors} \cite{HAN06}. In particular, in \cite{NAG03} 
some classical gradient methods have been reformulated as a linear combination of the singular vectors of $\mybf{A}$, 
in which the coefficients form involves the inverse of the singular values $\sigma_i$ multiplied for suitable factors that balance 
possible noise amplification effects due to small values of some $\sigma_i$. In this paper we want to go one step further, 
proving that a similar expansion holds true even if a scaling matrix multiplying the gradient is present. Moreover, 
we investigate the effect on the filter factors of some scaling matrices, arising from both the numerical optimization and 
the imaging framework. In particular, we show that the apparently negative oscillating 
behaviour of the filters provided by some scaled schemes results to be a surplus value for a more detailed reconstruction 
of the true solution. Finally, we extend our analysis to a particular class of projected gradient methods in the case of 
non-negative solutions.\\
The plan of the paper is the following: in Section 2 we introduce the optimization methods for the solution of \eqref{eq1} 
we decided to analyze, while in Section 3 the state-of-the-art on such algorithms as filtering methods is summarized, and 
the extensions to the scaled and projected cases are provided. Section 4 shows the behaviour of the 
filter factors for the considered regularization algorithms through some numerical experiments, while our conclusions are 
offered in Section 5.

\section{Gradient methods}\label{sec2}

A gradient method for the solution of \eqref{eq1} is an iterative algorithm whose $(k+1)$--th element is defined by
\begin{equation}\label{eq:SGiteration}
\mybf{x}_{k+1} = \mybf{x}_{k} -\alpha_k\mybf{M}_k\mybf{g}_k, 
\end{equation}
where:
\begin{itemize}
 \item $\mybf{g}_k = \nabla f(\mybf{x}_k) = \mybf{A}^T(\mybf{A}\mybf{x}_k - \mybf{b})$ is the gradient vector;
 \item $\mybf{M}_k$ is a symmetric and positive definite scaling matrix;
 \item $\alpha_k$ is the steplength parameter.
\end{itemize}
In the literature, a large variety of possibilities for choosing the steplength $\alpha_k$ and the scaling matrix $\mybf{M}_k$ 
has been proposed. Classical examples of steplength selection are the Steepest Descent (SD) \cite{CAU47} and the Minimal Gradient 
(MG) \cite{DAI03,ZHO06} methods, which minimize $f(\mybf{x}_k - \alpha\mybf{M}_k\mybf{g}_k)$ and 
$\|\nabla f(\mybf{x}_k - \alpha\mybf{M}_k\mybf{g}_k)\|_2$, respectively:
\begin{equation}\label{eq:aSD}
\alpha_k^{SD} = \underset{\alpha \in \mathbb{R}}{\text{argmin}} \ f(\mybf{x}_k - \alpha\mybf{M}_k\mybf{g}_k) = 
\frac{\mybf{g}_k^T\mybf{M}_k\mybf{g}_k}{\|\mybf{A} \mybf{M}_k\mybf{g}_k\|_2^2}.
\end{equation}
\begin{equation}\label{eq:aMG}
\alpha_k^{MG} = \underset{\alpha \in \mathbb{R}}{\text{argmin}} \ \|\nabla f(\mybf{x}_k - \alpha\mybf{M}_k\mybf{g}_k)\|_2 = 
\frac{\mybf{g}_k^T\mybf{M}_k \mybf{A}^T\mybf{A} \mybf{g}_k}{\|\mybf{A}^T\mybf{A}\mybf{M}_k\mybf{g}_k\|_2^2}.
\end{equation}
In order to accelerate the slow convergence exhibited in most cases by the standard formulas \eqref{eq:aSD} and \eqref{eq:aMG}, 
many other strategies for the steplength selection have been proposed, as the two Barzilai and Borwein (BB) rules \cite{BAR88} 
\begin{equation*}\label{eq:BB_non_scal}
\alpha_k^{BB1} = \frac{\mybf{s}_{k-1}^T\mybf{s}_{k-1}}{\mybf{s}_{k-1}^T\mybf{y}_{k-1}} \qquad ; \qquad
\alpha_k^{BB2} = \frac{\mybf{s}_{k-1}^T\mybf{y}_{k-1}}{\mybf{y}_{k-1}^T\mybf{y}_{k-1}},
\end{equation*}
where $\mybf{s}_{k-1} = \mybf{x}_{k} - \mybf{x}_{k-1}$ and $\mybf{y}_{k-1} = \mybf{g}_{k}-\mybf{g}_{k-1}$. These two schemes have been 
adapted to account for the scaling matrix \cite{BON09} as follows
\begin{equation*}\label{eq:BB}
\alpha_k^{BB1S} = \frac{\mybf{s}_{k-1}^T\mybf{M}_k^{-1}\mybf{M}_k^{-1}\mybf{s}_{k-1}}{\mybf{s}_{k-1}^T\mybf{M}_k^{-1}\mybf{y}_{k-1}} \qquad ; \qquad
\alpha_k^{BB2S} = \frac{\mybf{s}_{k-1}^T\mybf{M}_k\mybf{y}_{k-1}}{\mybf{y}_{k-1}^T\mybf{M}_k\mybf{M}_k\mybf{y}_{k-1}}.
\end{equation*}
Further accelerations have been proposed hereafter by alternating different steplength rules by means of an adaptively controlled 
switching criterion. Examples of such rules are the Adaptive Steepest Descent (ASD) method, the Adaptive Barzilai-Borwein (ABB) 
method \cite{ZHO06} and its generalizations $\text{ABB}_{\text{min1}}$ and $\text{ABB}_{\text{min2}}$ provided by Frassoldati et 
al. \cite{FRA08}.\\
As far as the scaling matrix concerns, the first example we considered is given by the well-known conjugate gradient algorithm for the 
least squares problem (CGLS). Indeed, this method can be expressed in the form \eqref{eq:SGiteration} by choosing the SD steplength 
and the scaling matrix
\begin{equation*} 
\mybf{M}_k = \mybf{I} - \frac{\mybf{s}_{k-1}\mybf{y}_{k-1}^T}{\mybf{y}_{k-1}^T\mybf{s}_{k-1}}. 
\end{equation*}
In our analysis, we also included two scaling matrices arising from the constrained optimization case. The first is the one provided by the iterative space reconstruction algorithm (ISRA), proposed by 
Daube-Whiterspoon and Muehllener \cite{DAU86} for reducing the computational cost of the expectation-maximization (EM) algorithm 
of Shepp and Vardi \cite{SHE82} in emission tomography. The explicit expression of the scaling matrix is
\begin{equation*}\label{ISRA}
\mybf{M}_k = \text{diag}\left(\frac{\mybf{x}_k}{\mybf{A}^T\mybf{A}\mybf{x}_k}\right),
\end{equation*}
where the quotient is intended in the Hadamard sense, and we obtain the positive definiteness by thresholding the diagonal elements 
in a prefixed interval $0 < L_{\min} < L_{\max}$.\\
The second scaling matrix is the one described in \cite{HAG09} by Hager, Mair and Zhang (HMZ), that exploits a gradient splitting 
strategy \cite{LAN01,LAN02} with a resulting scaling matrix given by 
\begin{equation}\label{Hager}
\mybf{M}_k = \text{diag}\left(\frac{\alpha_k^{CBB1}\mybf{x}_k}{\mybf{x}_k + \alpha_k^{CBB1}(\mybf{A}^T\mybf{A}\mybf{x}_k-
\mybf{A}^T\mybf{b})^+}\right),
\end{equation}
where $t^+ = \max\{0,t\}$ and $\alpha_k^{CBB1}$ is a cyclic version of the first BB steplength rule computed by reusing 
$\alpha_k^{BB1}$ for $p$ consecutive iterations \cite{DAI06}. A further thresholding step assures again the positive definiteness 
of the scaling matrix.\\
We remark that here we are {\em not} considering the ISRA and HMZ algorithms, but we only borrow the scaling
matrices defined in the algorithms themselves and use them in our scaled gradient scheme.\\
We point out that, while the convergence of the nonscaled gradient methods with the steplengths described before is guaranteed, 
the presence of a scaling matrix might require a thresholding of the steplength in a fixed range of positive values 
$[\alpha_{min},\alpha_{max}]$, followed by the introduction of a successive steplength reduction strategy. A classical example is the well-known Armijo rule \cite{BER99}: for given scalars $0 <  \gamma,\beta < 1$ and $z_k > 0$, the steplength $\alpha_k$ 
is set equal to $\beta^{q_k}z_k$, where $q_k$ is the first non-negative integer $q$ for which
\begin{equation*}\label{Armijo}
f(\mybf{x}_k)-f(\mybf{x}_k-\beta^qz_k\mybf{M}_k\nabla f(\mybf{x}_k)) \geq \gamma\beta^qz_k \nabla f(\mybf{x}_k)^T\mybf{M}_k\nabla f(\mybf{x}_k).
\end{equation*}

\section{Filtering effect}

When considering a linear inverse problem with noisy data
\begin{equation*} 
\mybf{b} = \mybf{A}\mybf{x}_{\text{true}} + \mybf{\eta}, 
\end{equation*}
it is well-known that the role played by a regularization method is to contrast the amplifying effect on the noise $\mybf{\eta}$ due to the small singular values of $\mybf{A}$. Indeed, if $\mybf{A} = \mybf{U}\mybf{\Sigma}\mybf{V}^T$ is the singular value decomposition of the matrix $\mybf{A}$, $\mybf{u}_i$, $\mybf{v}_i$ are the columns of $\mybf{U}$, $\mybf{V}$ and $\sigma_1 \geq \ldots \geq \sigma_n > 0$ are the diagonal elements of $\mybf{\Sigma}$, then the following relation for the solution $\mybf{x}^\dagger$ of \eqref{eq1} holds:
\begin{equation}\label{lincom}
\mybf{x}^\dagger = \sum_{i=1}^n \frac{\mybf{u}_i^T\mybf{b}}{\sigma_i}\mybf{v}_i = \mybf{x}_{\text{true}} + \sum_{i=1}^n \frac{\mybf{u}_i^T\mybf{\eta}}{\sigma_i}\mybf{v}_i.
\end{equation}
The classical recipe to contain the error propagation on the regularized solution $\mybf{x}_{\text{reg}}$ consists of adding some weights in the linear combination \eqref{lincom} that filter out the last components:
\begin{equation}\label{solreg}
\mybf{x}_{\text{reg}} = \sum_{i=1}^n \phi_i \frac{\mybf{u}_i^T\mybf{b}}{\sigma_i}\mybf{v}_i.
\end{equation}
The truncated singular value decomposition (TSVD) and the Tikhonov method are examples of this approach, and the corresponding regularized solutions can be written in the form \eqref{solreg} with, respectively, 
\begin{equation*} 
\phi_i = \begin{cases} 1 & \text{if } \ i \leq r \\ 0 & \text{if } \ i > r \end{cases} \qquad ; \qquad \phi_i = \frac{\sigma_i^2}{\sigma_i^2 + \lambda}, \end{equation*}
where $r \in \{1,\ldots,n\}$ and $\lambda > 0$.

\subsection{The nonscaled case}

\begin{figure}[ht]
\begin{center}
 \begin{tabular}{ccc}
\includegraphics[scale = 0.27]{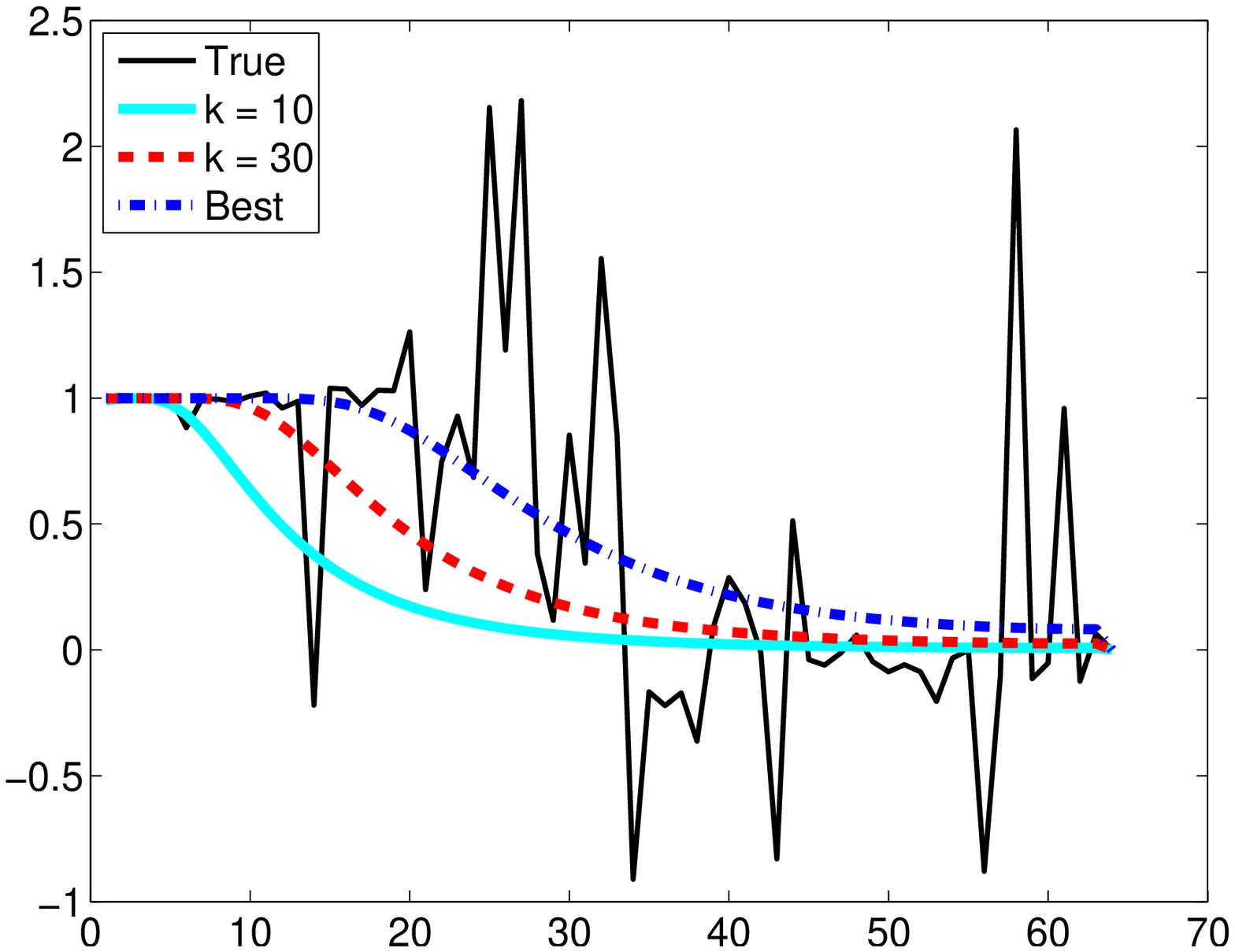}&
\includegraphics[scale = 0.27]{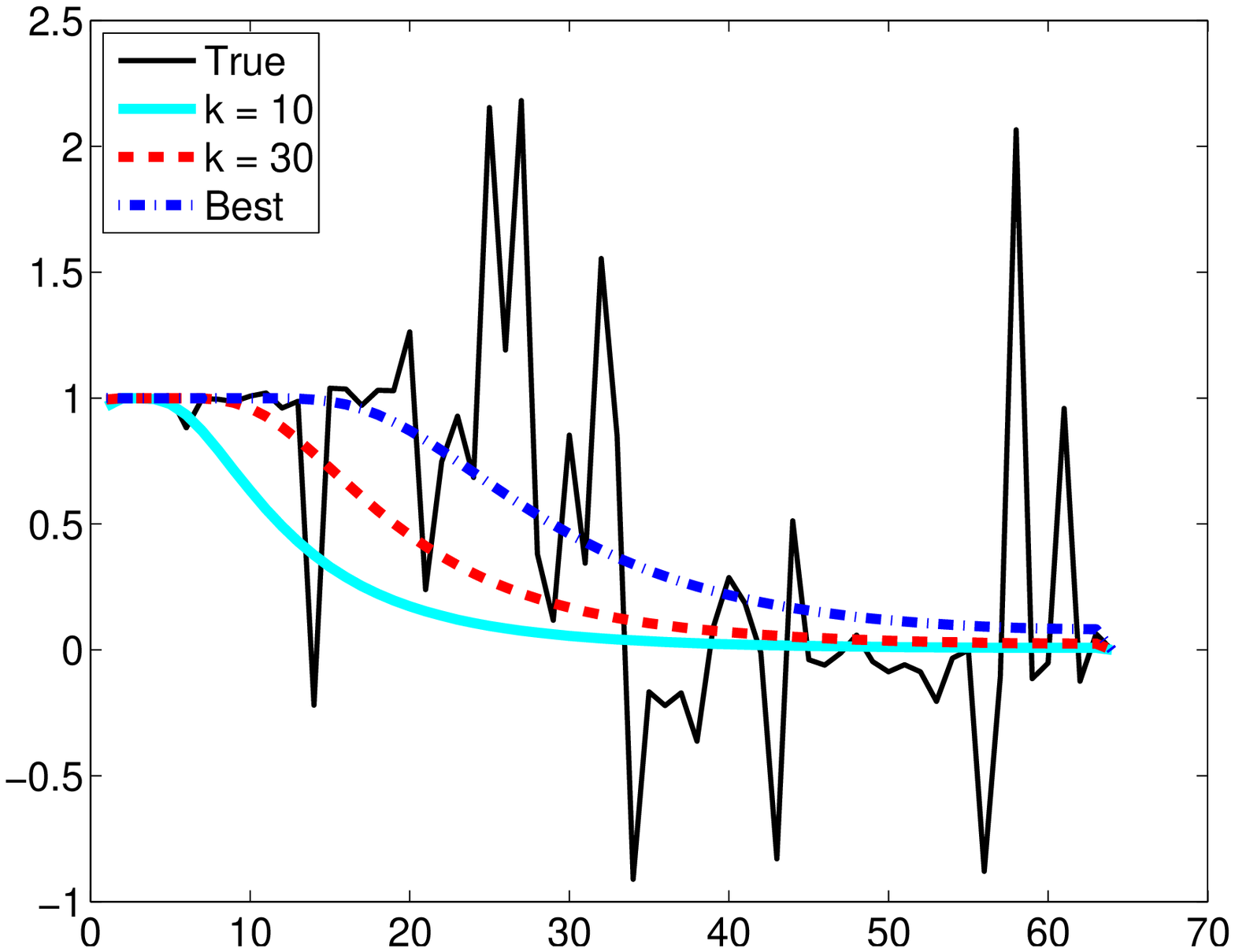}&
\includegraphics[scale = 0.27]{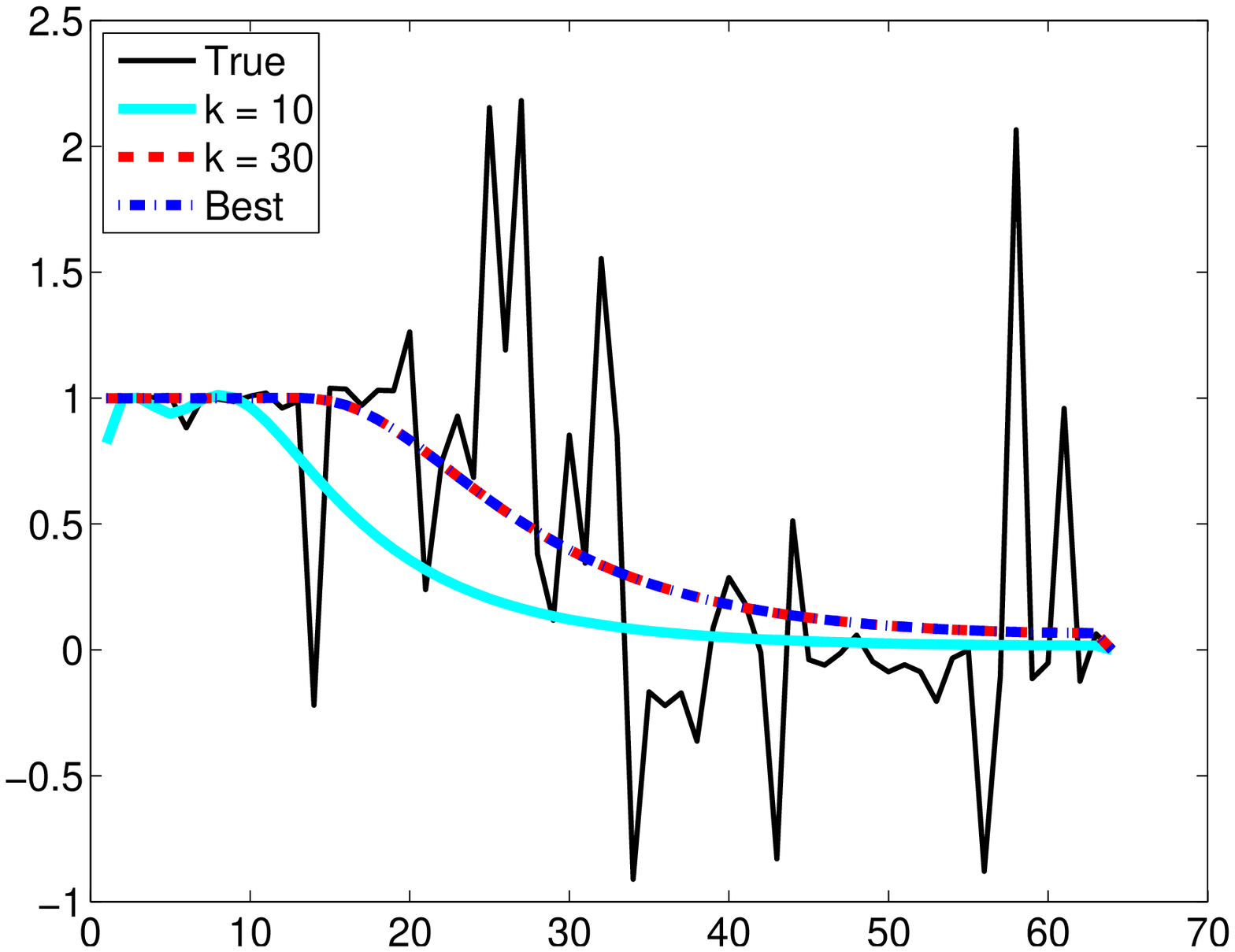}\\
MG & SD & BB1\\
\includegraphics[scale = 0.27]{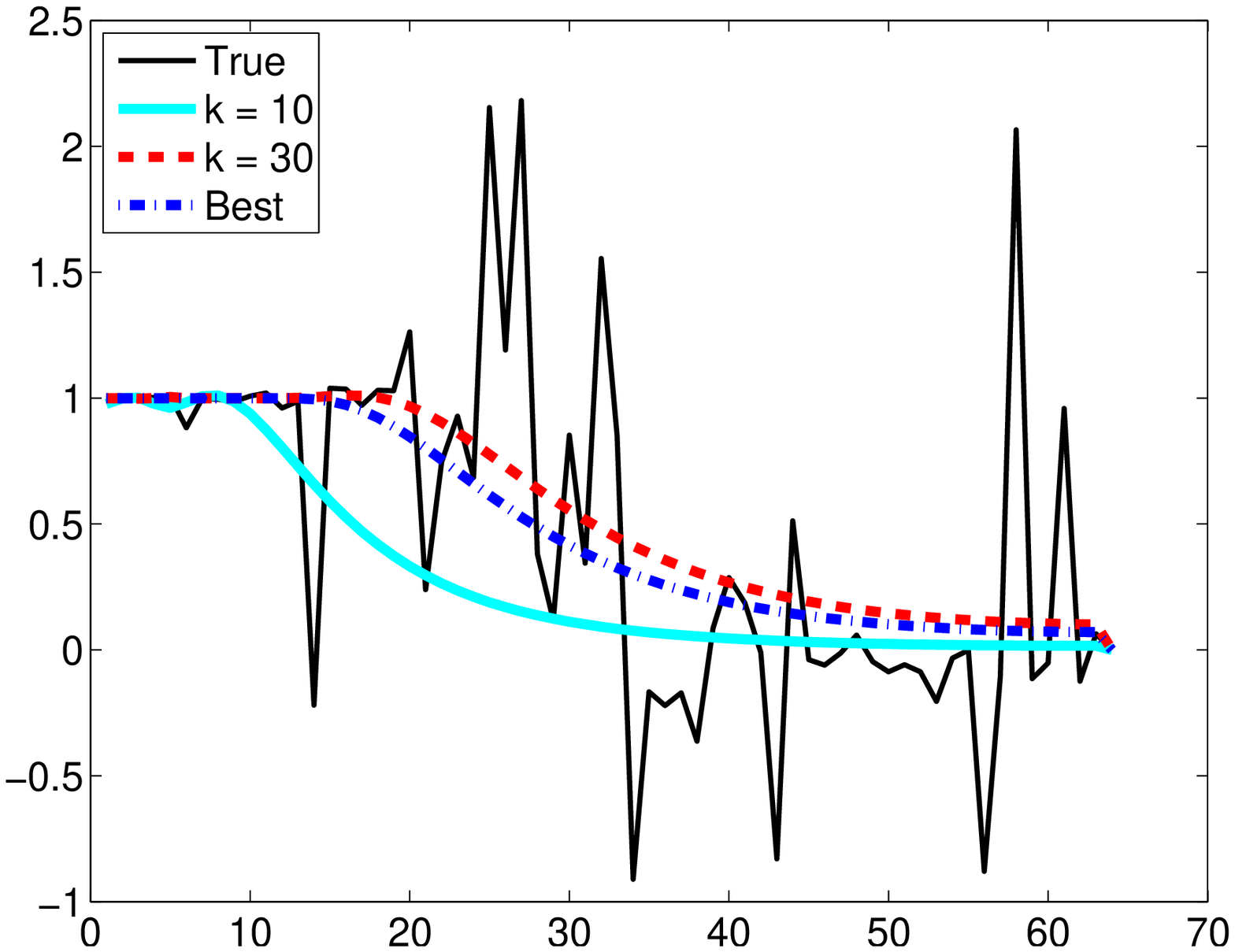}&
\includegraphics[scale = 0.27]{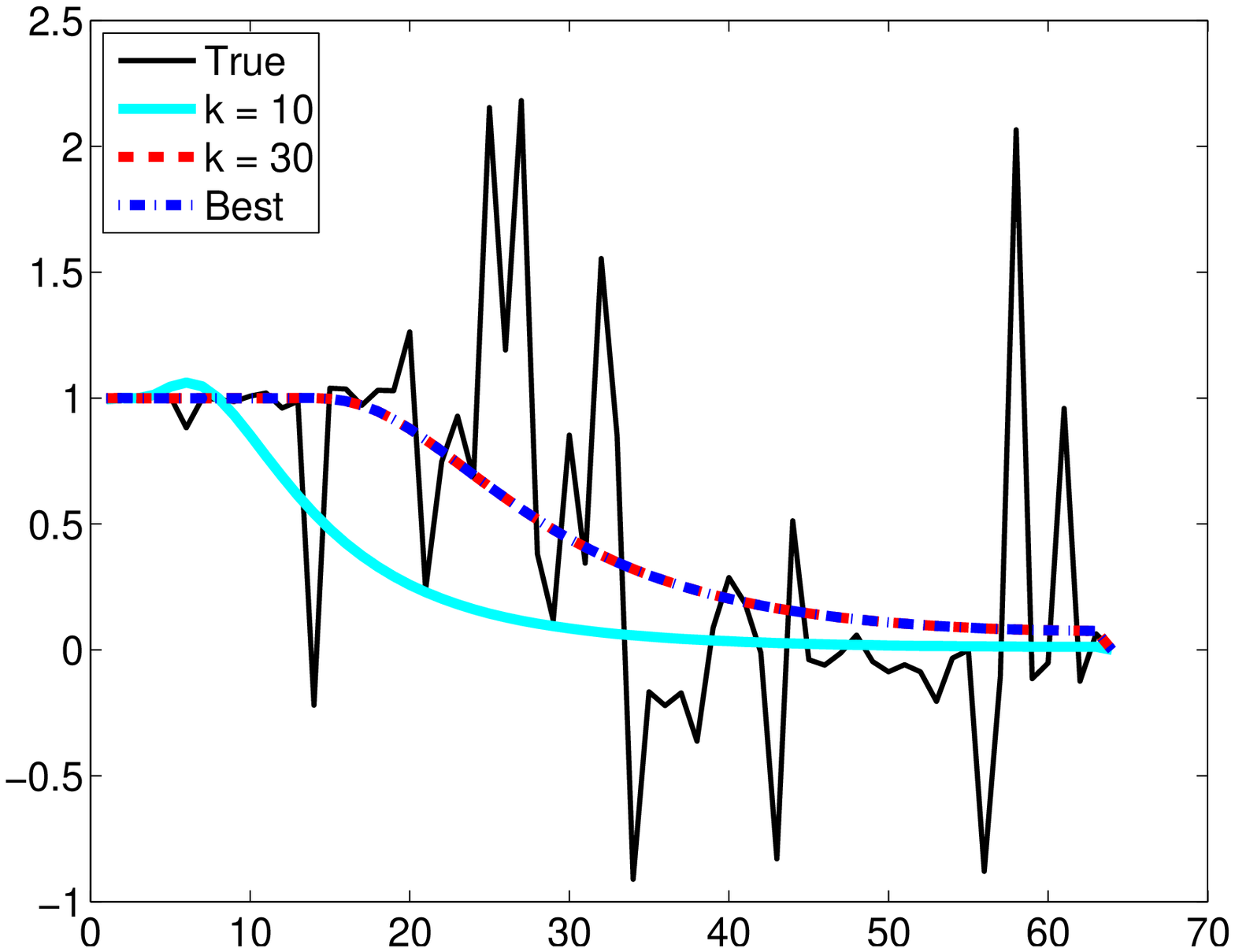}&
\includegraphics[scale = 0.27]{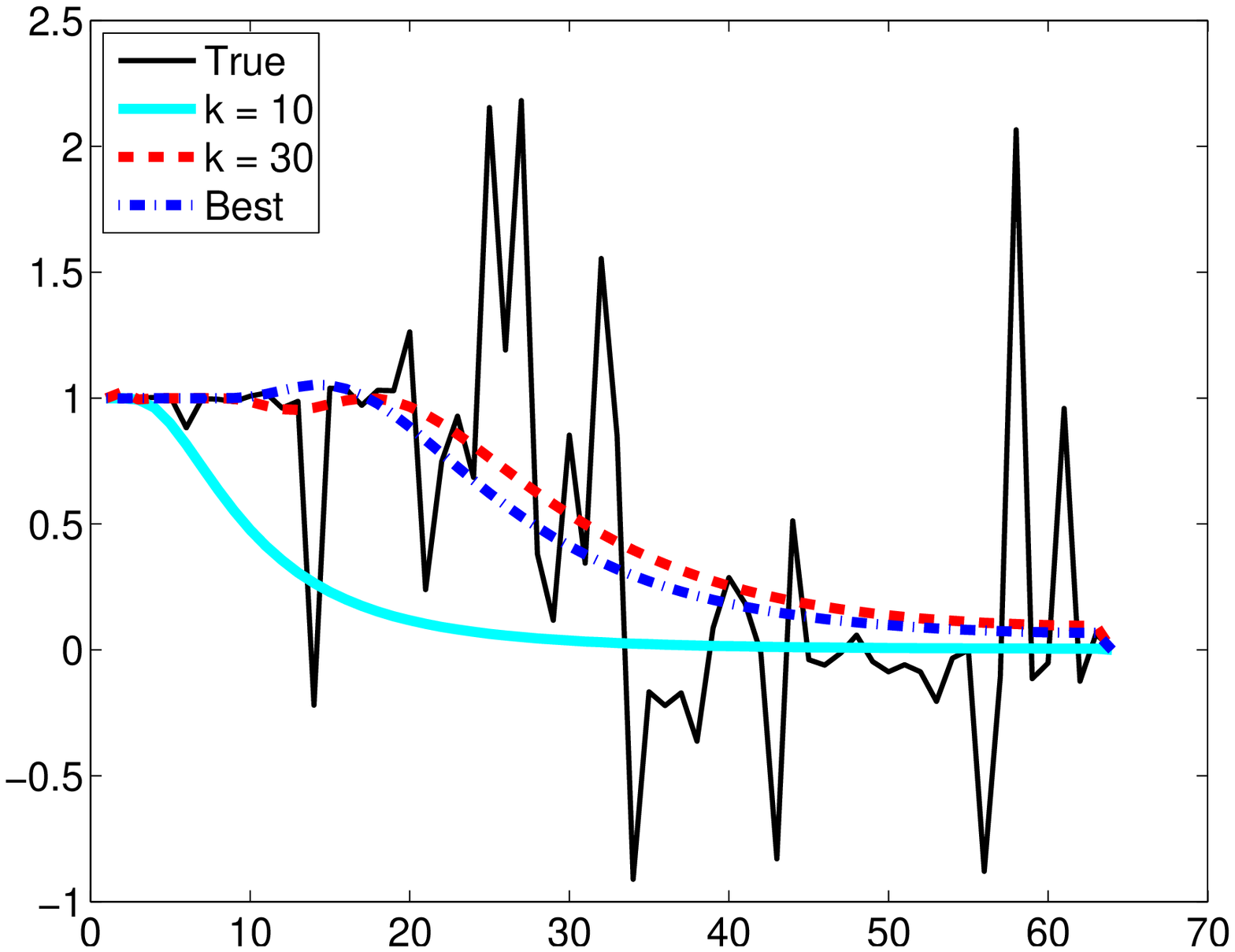}\\
BB2  & ABB & ABB$_{\text{min1}}$
 \end{tabular}
\caption{Comparison of the $\mybf{x}_{\text{true}}$ filter factors (black thin) for the \texttt{heat} test problem with the ones generated by gradient methods with different steplength rules at iteration 10 (cyan thick), 30 (red dashed) and at the iteration corresponding to the minimum error (blue dashdot).}\label{fig:FiltriHeatNoScal}
\end{center}
\end{figure}

Besides the classical ``direct'' approaches as TSVD and Tikhonov, also iterative regularization strategies can be thought by means of their filtering effect. Indeed, expression \eqref{solreg} can be generalized to the regularized solution $\mybf{x}_{k+1}$ provided by any gradient method by writing 
\begin{equation*}\label{solregk}
\mybf{x}_{k+1} = \sum_{i=1}^n \phi_i^{k+1} \frac{\mybf{u}_i^T\mybf{b}}{\sigma_i}\mybf{v}_i.
\end{equation*}
where the filter factors $\phi_i^{k+1}$ are automatically defined during the iterative procedure. The easiest example of gradient method is the Landweber algorithm, which generates a sequence $\{\mybf{x}_k\}$ by
\begin{equation}\label{eq:Landweber}
\mybf{x}_{k+1} = \mybf{x}_{k} -\alpha\mybf{g}_k = \mybf{x}_k -\alpha\mybf{A}^T(\mybf{A}\mybf{x}_k - \mybf{b})\,,
\end{equation}
where the steplength $\alpha$ is fixed during the iterations and must satisfy $0 < \alpha < 2/\sigma_1^2$ in order
to guarantee the convergence. The iteration \eqref{eq:Landweber} can be rewritten, by choosing the initial guess 
$\mybf{x}_0 = \mybf{0}$, as
\begin{equation*}
\mybf{x}_{k+1} = \alpha\sum_{\ell = 0}^k(\mybf{I} - \alpha\mybf{A}^T\mybf{A})^\ell\mybf{A}^T\mybf{b}
\end{equation*}
or, equivalently, by means of the SVD as
\begin{equation*}
\mybf{x}_{k+1} = \sum_{i = 1}^n\left(1 - (1 -\alpha\sigma_i^2)^{k+1}\right)\frac{\mybf{u}_i^T\mybf{b}}{\sigma_i}\mybf{v}_i\,,
\end{equation*}
thus obtaining an expression for the filter factors given by
\begin{equation*}\label{eq:FFLandweber}
\phi_i^{k+1} = 1 - (1 - \alpha\sigma_i^2)^{k+1}\quad k=0,1,\ldots
\end{equation*}
The previous equation suggests an extension to the case of the steplength $\alpha_k$ varying at each iteration, in the case of nonscaled gradient. If $\mybf{M}_k = \mybf{I}$, indeed, the filter factors of a general gradient method can be written as
\begin{equation}\label{eq:FFgradient-type_methods}
\phi_i^{k+1} = 1 - \prod_{\ell = 0}^k \bigl(1 - \alpha_\ell\sigma_i^2\bigr)\quad k=0,1,\ldots
\end{equation}
Since the sequence generated by the gradient method converges to the least-squares solution $\mybf{x}^\dagger$, the filter factors 
$\phi_i^{k+1}$ will tend to 1 as $k$ increases. Moreover, from the form of the products in \eqref{eq:FFgradient-type_methods} we 
can observe that the convergence rate of the sequence $\{\phi_i^{k+1}\}$ to its limit value depends on the steplength values 
$\alpha_\ell$. In Figure \ref{fig:FiltriHeatNoScal}, the filter factors obtained in a numerical experiment with different steplength 
rules have been reported as functions of the singular values index $i$. In particular, we considered the \texttt{heat} dataset 
from Hansen's \textit{Regularization Tools} \cite{HAN94} and we plotted the filter factors generated at iteration 10, 30 and at 
the iteration corresponding to the minimum error for gradient methods with the MG, SD, BB1, BB2, ABB and ABB$_{\text{min1}}$ 
steplength rules (in this very last case, we used the revised version of the ABB$_{\text{min1}}$ scheme 
proposed in \cite{PRA12}). In all cases the filter factors start from zero (since we chose $\mybf{x}_0 = \mybf{0}$ as initial point) and 
converge to 1, with a faster convergence clearly visible for the factors corresponding to the largest singular values.
Moreover, the regular and smooth filter trends achieved in all the figures show that any steplength selection rule is not 
able to generate filter factors that behave as erratically as the true ones given by 
$\sigma_i(\mybf{v}_i^T\mybf{x}_{\text{true}})/(\mybf{u}_i^T\mybf{b})$.

\subsection{Adding a scaling matrix}

If we assume $\mybf{x}_0 = \mybf{0}$, then the general expression \eqref{eq:SGiteration} of a gradient method can be easily rewritten in the form 
\begin{equation*}
\mybf{x}_{k+1} = P_k(\mybf{A}^T\mybf{A})\mybf{A}^T\mybf{b},
\end{equation*}
where the polynomial $P_k$ acts on a $n \times n$ matrix $\mybf{\Omega}$ as
\begin{equation}\label{matpol}
P_k(\mybf{\Omega}) = P_{k-1}(\mybf{\Omega})+\alpha_k\mybf{M}_k(\mybf{I} - \mybf{\Omega} P_{k-1}(\mybf{\Omega})), \qquad P_{-1}(\mybf{\Omega}) = \mybf{0}.
\end{equation}
By writing $\mybf{x}_{k+1}$ in terms of the basis formed by $\mybf{v}_i$ and setting $\mybf{Q}_k = \mybf{V}^T P_k(\mybf{A}^T\mybf{A})\mybf{V}$, the following equations hold:
\begin{eqnarray}
\label{eq:xfilt_scal}
\nonumber
\mybf{x}_{k+1} & = & \mybf{V}\mybf{V}^T P_k(\mybf{A}^T\mybf{A})\mybf{V}\mybf{V}^T \mybf{A}^T\mybf{b} = \nonumber \\
& = &\sum_{i=1}^n (\mybf{Q}_k\mybf{V}^T\mybf{A}^T\mybf{b})_i\mybf{v}_i = \sum_{i=1}^n (\mybf{Q}_k\mybf{\Sigma}^T\mybf{U}^T\mybf{b})_i\mybf{v}_i =\nonumber \\ 
& = & \sum_{i=1}^n \left( \sum_{j=1}^n (\mybf{Q}_k)_{ij}\sigma_j \mybf{u}_j^T\mybf{b}\right) \mybf{v}_i =\nonumber\\
& = & \sum_{i=1}^n\sigma_i^2 \left( \frac{\sum_{j=1}^n (\mybf{Q}_k)_{ij}\sigma_j\mybf{u}_j^T\mybf{b}}{\sigma_i\mybf{u}_i^T\mybf{b}}\right) \frac{\mybf{u}_i^T\mybf{b}}{\sigma_i}\mybf{v}_i.
\end{eqnarray}
Comparing \eqref{solreg} and \eqref{eq:xfilt_scal}, it is easy to see that the filter factors for the $(k+1)$-th iteration of a scaled gradient method have the expression below:
\begin{equation}\label{filtscal}
\phi_i^{k+1}= \sigma_i^2 \left( \frac{\sum_{j=1}^n(\mybf{Q}_k)_{ij}\sigma_j \mybf{u}_j^T\mybf{b}}{\sigma_i\mybf{u}_i^T\mybf{b}}\right)\,,\qquad i=1,\ldots,n.
\end{equation}
Since, when $\mybf{M}_k$ is the identity matrix, $\mybf{Q}_k$ can be written as
\begin{align*}
\mybf{Q}_k = & \mybf{V}^TP_k(\mybf{A}^T\mybf{A})\mybf{V} = \mybf{V}^TP_k(\mybf{V}\mybf{\Sigma}^T\mybf{\Sigma}\mybf{V}^T)\mybf{V} =\\
= & \mybf{V}^T\left[\mybf{V} P_k(\mybf{\Sigma}^T\mybf{\Sigma})\mybf{V}^T\right]\mybf{V} = P_k(\mybf{\Sigma}^T\mybf{\Sigma})\,,
\end{align*}
it follows that the general expression \eqref{filtscal} for the nonscaled case becomes, according to
\cite{NAG03},
\begin{equation}
\label{filtnonscal}
\phi_i^{k+1} = \sigma_i^2P_{k+1}(\sigma_i^2)\,, \qquad i=1,\dots,n,
\end{equation}
where equation \eqref{matpol} reduces to a usual scalar polynomial
\begin{equation*}
P_k(\xi) = P_{k-1}(\xi)+\alpha_k(1 - \xi P_{k-1}(\xi)), \qquad P_{-1}(\xi) = 0.
\end{equation*}
From equation \eqref{filtscal} we can see the effect of the scaling matrix on the filters expression: the presence of 
$\mybf{M}_k$ makes any factor $\phi_i^{k+1}$ related to all the singular values $\sigma_1,\ldots,\sigma_n$ and all 
the singular vectors $\mybf{u}_1,\ldots,\mybf{u}_n$, while for the nonscaled case each value depends only of the $i$-th singular 
value $\sigma_i$. For this reason, the filter factors expression \eqref{filtnonscal}, for the nonscaled methods, cannot be 
directly applied to the scaled case. In Section \ref{numexp}, we will show the positive effect 
of such more complicated dependence in reconstructing the actual values of the true solution filter factors.

\subsection{Non-negative solutions}\label{secproj}

When the unknown $\mybf{x}$ refers to a physical quantity, it is quite common to require its entries to be non-negative and, 
consequently, the problem \eqref{eq1} becomes
$$
\min_{\mybf{x} \geq 0} f(\mybf{x}) \equiv \frac{1}{2} \| \mybf{A} \mybf{x} - \mybf{b} \|_2^2.
$$
The modified general form of a gradient method accounting for this constraint on $\mybf{x}$ is the following \cite{BER99}:
\begin{equation}\label{gradproj2}
\mybf{x}_{k+1} = \mybf{x}_{k} + \lambda_k\left(\mathbb{P}^+_{\mybf{M}_k^{-1}}(\mybf{x}_{k} -\alpha_k\mybf{M}_k\mybf{g}_k)-\mybf{x}_{k}\right), \quad k=0,1,\ldots,
\end{equation}
where $\lambda_k$ is a linesearch parameter (typically chosen by means of an Armijo rule along the feasible direction) and $\mathbb{P}^+_{\mybf{M}}(\mybf{z})$ denotes the projection of the point $\mybf{z}$ on the non-negative orthant with respect to the norm induced by the matrix $\mybf{M}$, i.e.
\begin{equation*}
\mathbb{P}^+_{\mybf{M}}(\mybf{z}) = \underset{\mybf{x} \geq 0}{\text{argmin}} (\mybf{x}-\mybf{z})^T \mybf{M} (\mybf{x}-\mybf{z}).
\end{equation*}
For a general scaled gradient method, therefore, the introduction of a constraint on $\mybf{x}$ requires the solution of a quadratic program at each iteration. Here we consider the case in which the projection on the non-negative orthant can be expressed by means of a left multiplication for a diagonal matrix $\mybf{D}_k$, where
\begin{equation*}
(\mybf{D}_k)_{ii} = \begin{cases} 1 & \text{if } (\overline{\mybf{x}}_{k+1})_i \geq 0 \\ 0 & \text{if } (\overline{\mybf{x}}_{k+1})_i < 0 \end{cases}
\end{equation*}
and
\begin{equation*}
\overline{\mybf{x}}_{k+1} = \mybf{x}_{k} -\alpha_k\mybf{M}_k\mybf{g}_k.
\end{equation*}
The resulting projected gradient method, therefore, assumes the form
\begin{equation}\label{gradprojD}
\mybf{x}_{k+1} = \mybf{x}_{k} + \lambda_k(\mybf{D}_k\overline{\mybf{x}}_{k+1}-\mybf{x}_{k}), \quad k=0,1,\ldots
\end{equation}
The easiest (and mostly used) case of scaling matrix leading to an expression \eqref{gradprojD} for the $(k+1)$--th iteration is a diagonal $\mybf{M}_k$. For a projected gradient method in the form \eqref{gradprojD}, the expression of the filter factors is given again by \eqref{filtscal}, in which the polynomial $P_k$ defining the matrix $\mybf{Q}_{k}$ is generalized accounting for the presence of the projection matrix:
\begin{equation*}\label{matpolproj2}
P_k(\mybf{\Omega}) = P_{k-1}(\mybf{\Omega})+\lambda_k\left[\alpha_k\mybf{D}_k\mybf{M}_k - (\mybf{I} - \mybf{D}_k + \alpha_k\mybf{D}_k\mybf{M}_k\mybf{\Omega}) P_{k-1}(\mybf{\Omega}))\right], \quad P_{-1}(\mybf{\Omega}) = \mybf{0}.
\end{equation*}
A similar expression holds also for the specific case of constant linesearch parameter fixed at unity, $\lambda_k=1$, $k=0,1,\ldots$ (in this case, the convergence of the scheme is guaranteed by an Armijo rule along the projection arc \cite{BER99}). Under this assumption, the gradient iteration becomes
\begin{equation}\label{gradproj2D}
\mybf{x}_{k+1} = \mybf{D}_k\overline{\mybf{x}}_{k+1}, \quad k=0,1,\ldots,
\end{equation}
and the corresponding filter factors can be written again in the form \eqref{filtscal} with 
\begin{equation*}\label{matpolproj}
P_k(\mybf{\Omega}) = \mybf{D}_k\left(P_{k-1}(\mybf{\Omega})+\alpha_k\mybf{M}_k(\mybf{I} - \mybf{\Omega} P_{k-1}(\mybf{\Omega}))\right), \qquad P_{-1}(\mybf{\Omega}) = \mybf{0}.
\end{equation*}
The advantage of this simplified case is that expression \eqref{gradproj2D}, while being still valid for a diagonal $\mybf{M}_k$, can be extended also to a wider class of scalings, as the {\em diagonal matrices with respect to} $I^+(\mybf{x}_{k})$ \cite{BER99}, where
\begin{equation*}
I^+(\mybf{x}_{k}) = \left\{ i\in\{1,\ldots,n\} \ \big| \ (\mybf{x}_{k})_i=0,\frac{\partial f(\mybf{x}_{k})}{\partial x_i} > 0 \right\}.
\end{equation*}

\section{Numerical experiments}\label{numexp}

We show now the behaviour of the filter factors for the solutions obtained with the scaling matrices reported in Section \ref{sec2}, in the cases of 
both unconstrained and constrained problems. The evaluation of the results will be carried out on two different examples taken 
from Hansen's \textit{Regularization Tools} \cite{HAN94}. The first test is the \texttt{heat} 1-dimensional dataset already 
mentioned in the previous section. For this test problem we set up $n=64$ and assumed Gaussian noise with zero mean, variance 
equals to 1 and scaled so that $\|\mybf{\eta}\|_2 / \|\mybf{A}\mybf{x}_{\text{true}}\|_2 = 0.01$. In Table \ref{Tab1} we report 
the best reconstruction errors $\| \mybf{x}_k - \mybf{x}_{\text{true}} \|_2 / \| \mybf{x}_{\text{true}} \|_2$ for the six 
steplength rules considered in Figure \ref{fig:FiltriHeatNoScal} when $\mybf{D}_k = \mybf{M}_k = \mybf{I}$ for all $k$. 
From the results obtained, we can see that at the best iteration all the steplengths provide similar filter factors and, 
consequently, similar approximated solution. 

\begin{table}[ht]
\caption{Numbers of iterations required by nonscaled gradient methods with different steplengths and minimum error reached 
for solving \texttt{heat}.}
\label{Tab1}
\begin{center}
\begin{tabular}{l|cc}
\hline\noalign{\smallskip}
Method  						& Iter & Min err  \\
\noalign{\smallskip}\hline\noalign{\smallskip}
MG    							&   94 &    0.049 \\
SD    							&   97 &    0.049 \\
BB1   							&   30 &    0.049 \\
BB2   							&   29 &    0.049 \\
ABB   							&   30 &    0.048 \\
ABB$_\text{min1}$   &   29 &    0.048 \\
\noalign{\smallskip}\hline
\end{tabular}
\end{center}
\end{table}

\begin{figure}[ht]
\begin{center}
 \begin{tabular}{cc}
\includegraphics[scale = 0.30]{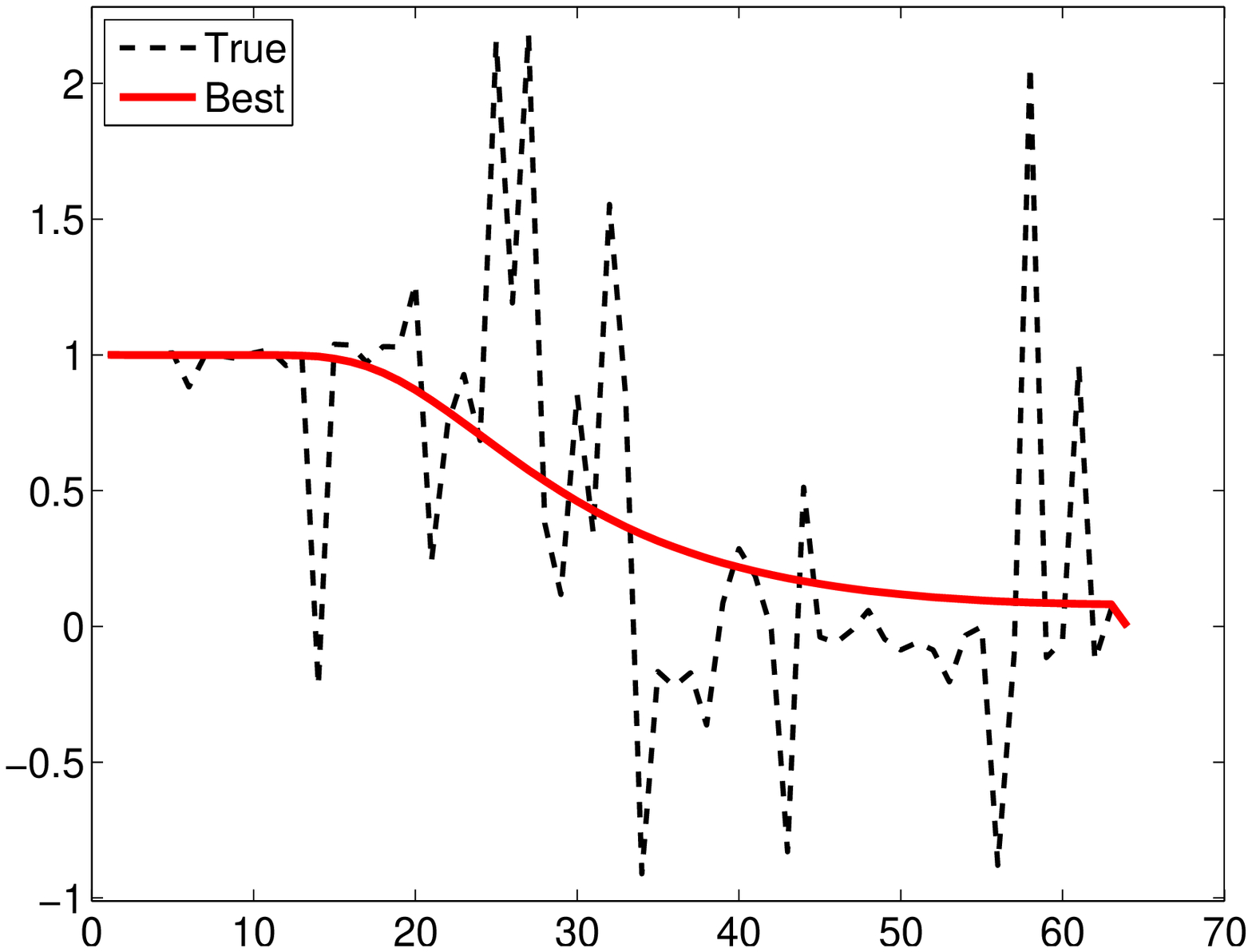}&
\includegraphics[scale = 0.30]{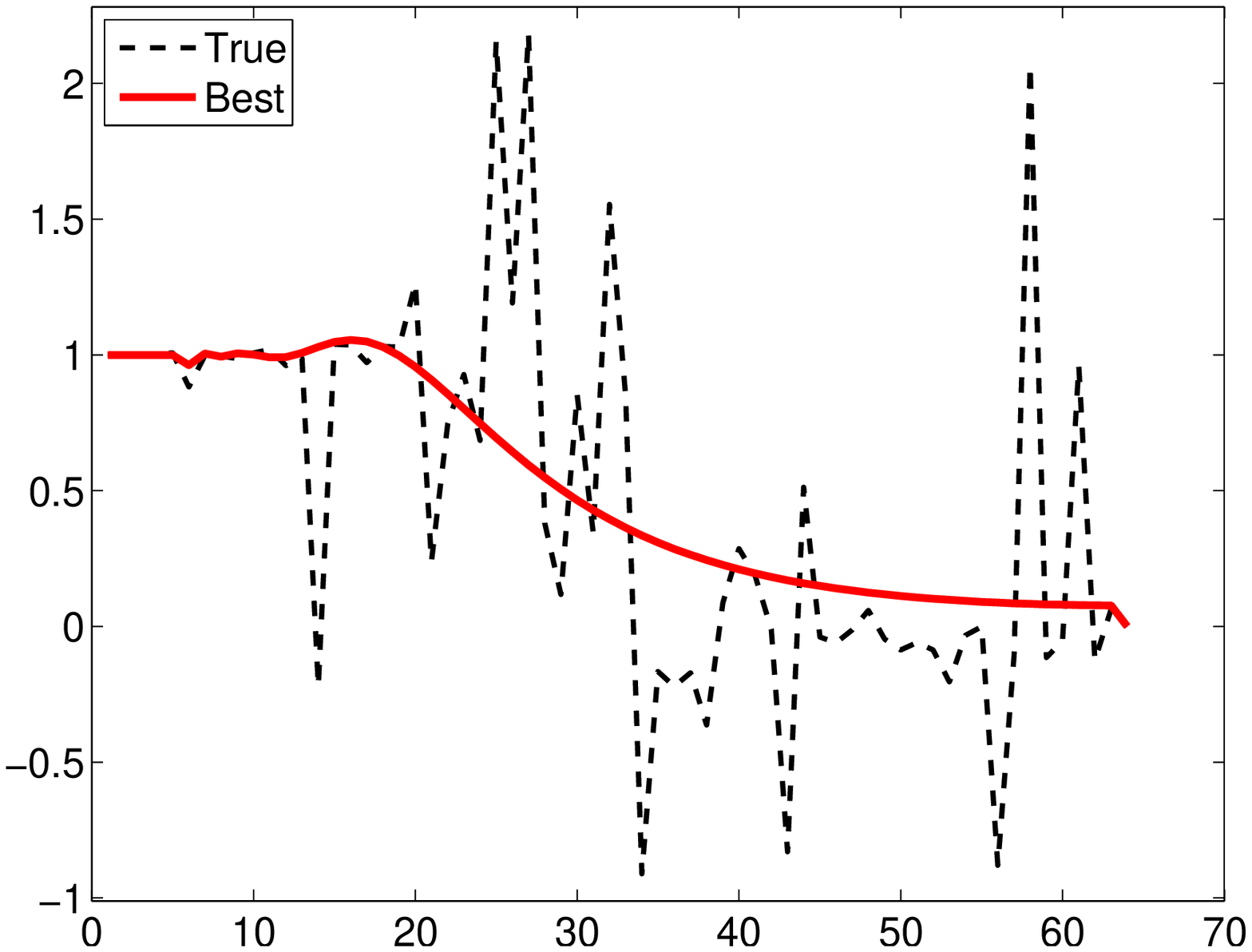}\\
SD & CGLS \\
\includegraphics[scale = 0.30]{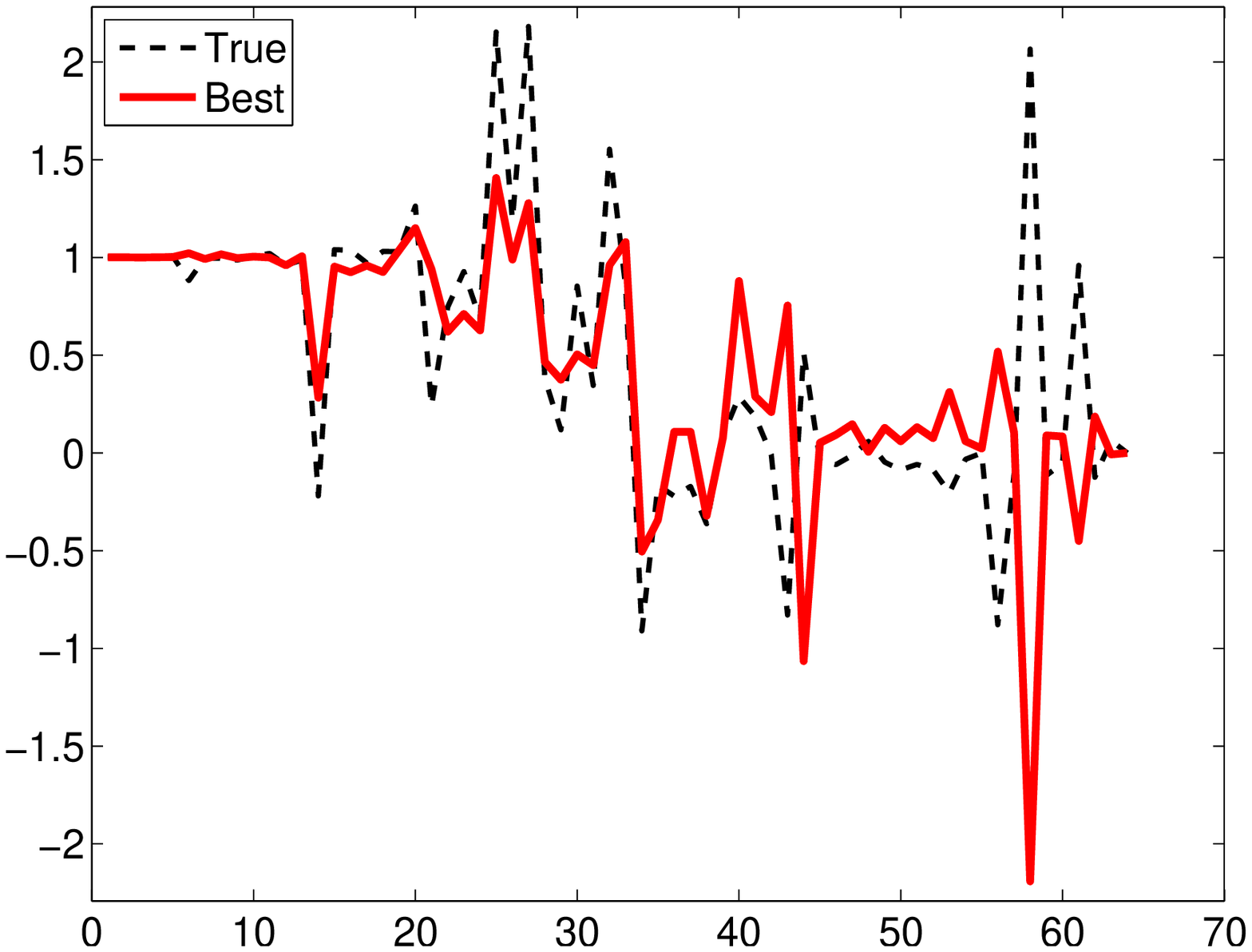}&
\includegraphics[scale = 0.30]{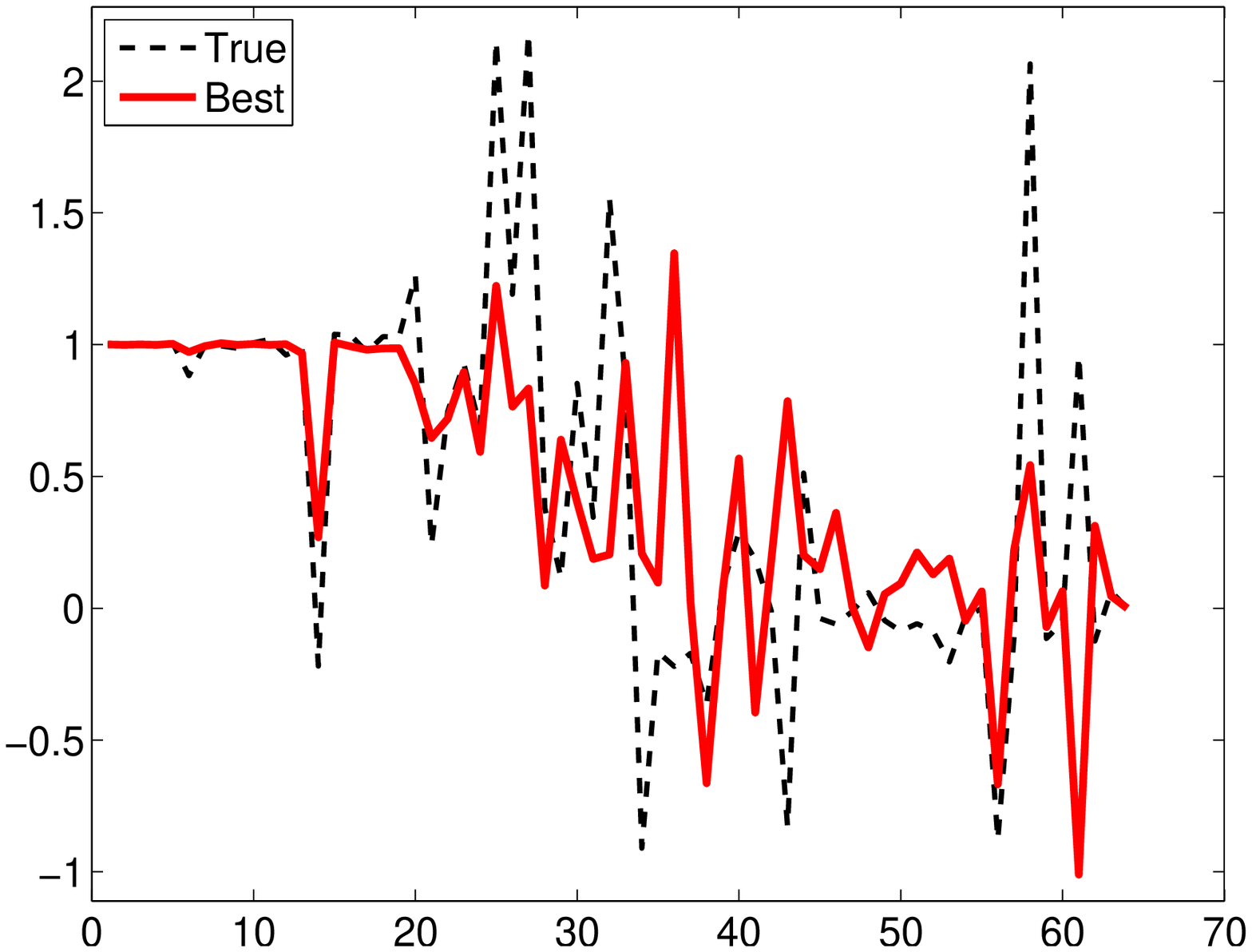}\\
ISRA & HMZ \\
 \end{tabular}
\caption{Comparison of the $\mybf{x}_{\text{true}}$ filter factors (black dashed) for the \texttt{heat} test problem 
with the ones generated by gradient methods with different scaling matrices at the iteration corresponding to the minimum 
error (red solid). The nonscaled filters are also reported (SD).}
\label{fig:FiltriHeatScal}
\end{center}
\end{figure}

\noindent For this reason, in order to evaluate the performance of the scaled and/or projected versions, in the following we 
will adopt the SD rule for the steplength selection. We consider first the unconstrained case and we analyze the behaviour of the different scalings 
on the restored solutions. The best reconstruction errors obtained by ISRA, HMZ and CGLS matrices are provided in Table \ref{Tab2}, 
while the corresponding filter factors are shown in Figure \ref{fig:FiltriHeatScal}. For the HMZ case we chose the value $p=4$ 
for the cyclic version of the BB1 steplength rule defining $\mybf{M}_k$ - see equation \eqref{Hager}. Moreover, the diagonal 
elements of the scaling matrices provided by the ISRA and HMZ approaches have been thresholded in the range $[10^{-3},10^{8}]$.

\begin{table}[ht]
\caption{Optimal numbers of iterations required by gradient methods with SD steplength and different scaling matrices 
for solving \texttt{heat}, with the corresponding minimum error. The suffix `\_P' denotes the projected algorithms.}
\label{Tab2}
\begin{center}
\begin{tabular}{l|cc}
\hline\noalign{\smallskip}
Method  & Iter & Min err \\
\noalign{\smallskip}\hline\noalign{\smallskip}
SD      &   97 &   0.049 \\
CGLS    &   11 &   0.047 \\
ISRA    &   25 &   0.037 \\
HMZ     &   33 &   0.044 \\
\noalign{\smallskip}\hline\noalign{\smallskip}
SD\_P   &  116 &   0.037 \\
ISRA\_P &   14 &   0.034 \\
HMZ\_P  &   66 &   0.037 \\
\noalign{\smallskip}\hline
\end{tabular}
\end{center}
\end{table}

\begin{figure}[ht]
\begin{center}
\begin{tabular}{ccc}
\includegraphics[scale = 0.27]{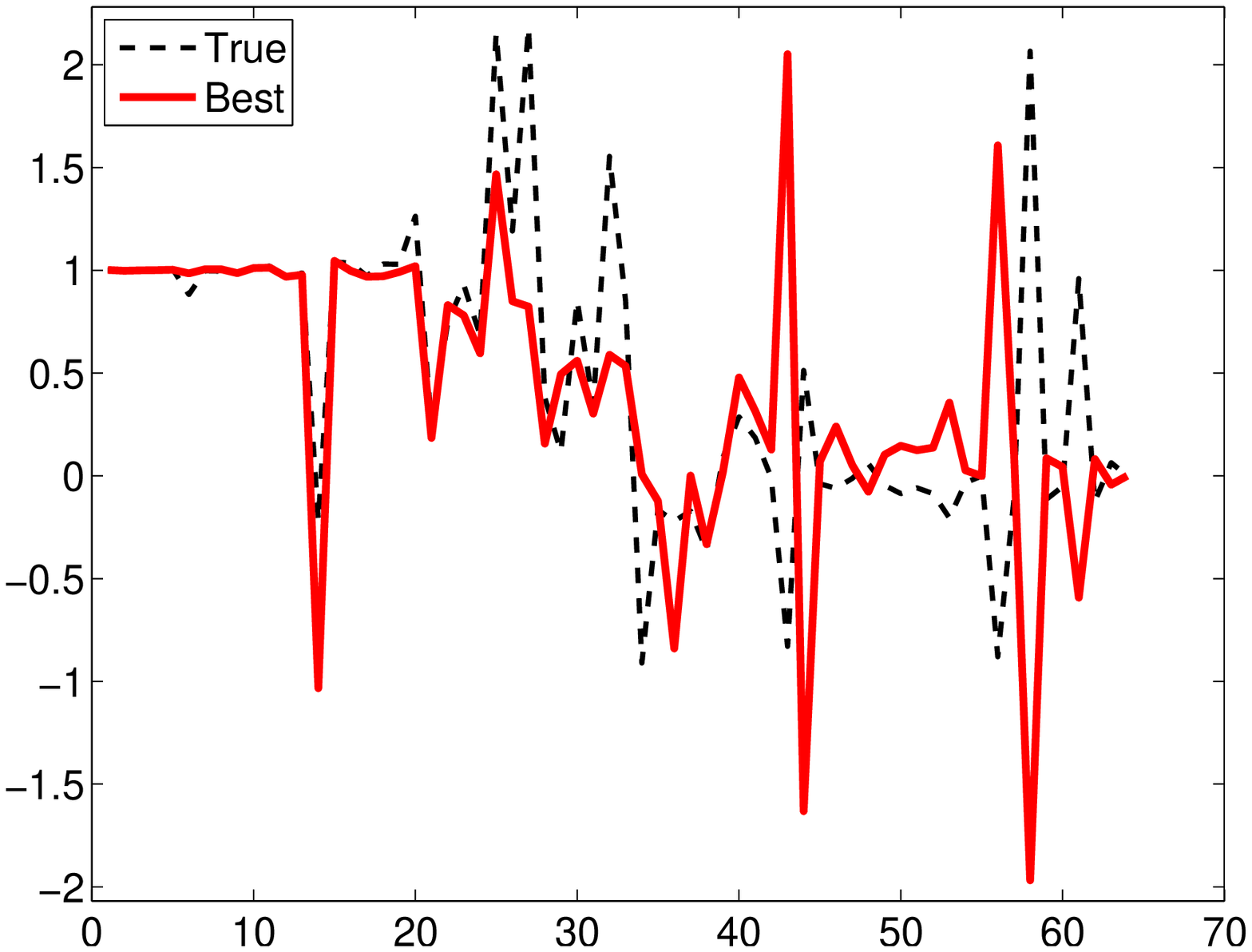}&
\includegraphics[scale = 0.27]{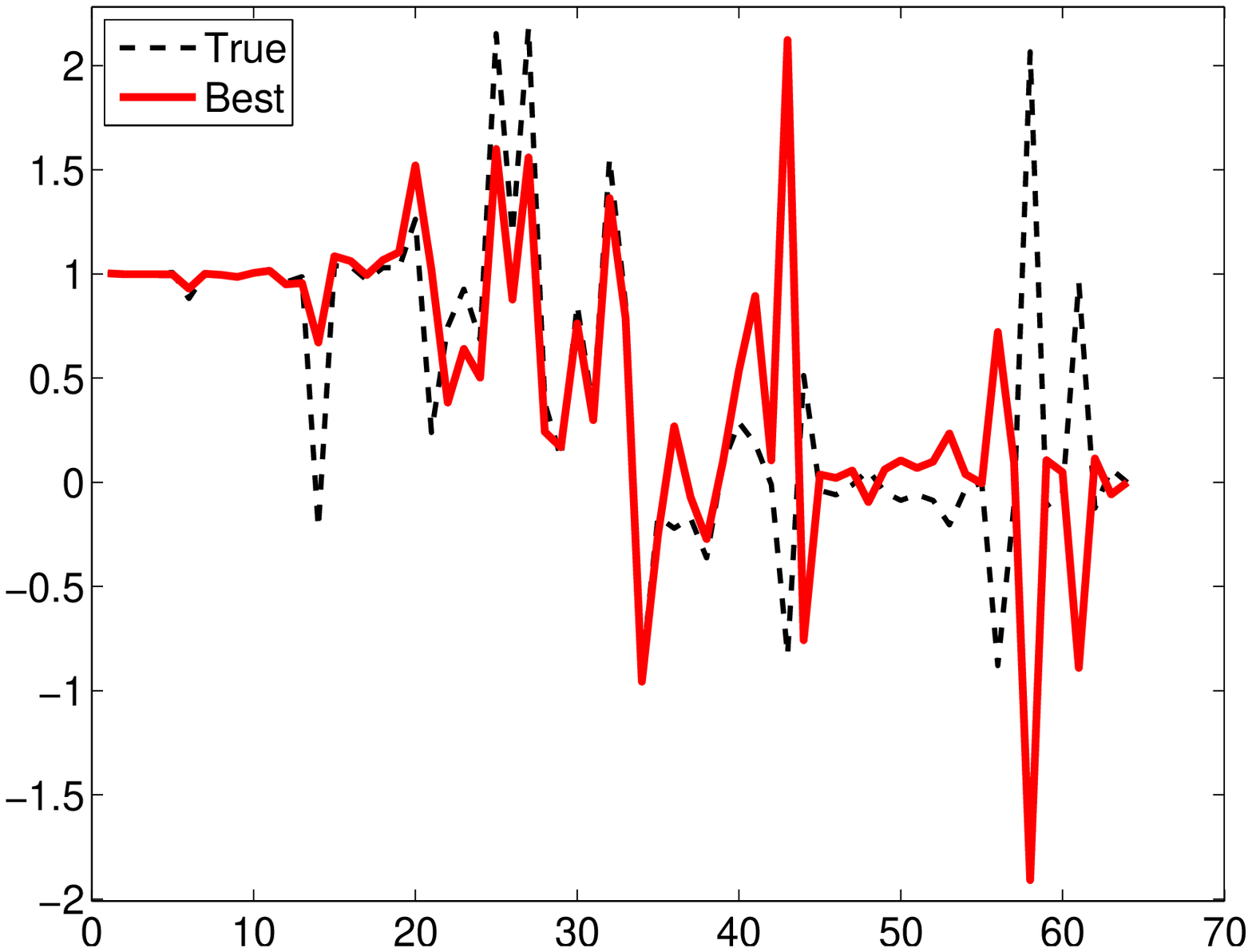}&
\includegraphics[scale = 0.27]{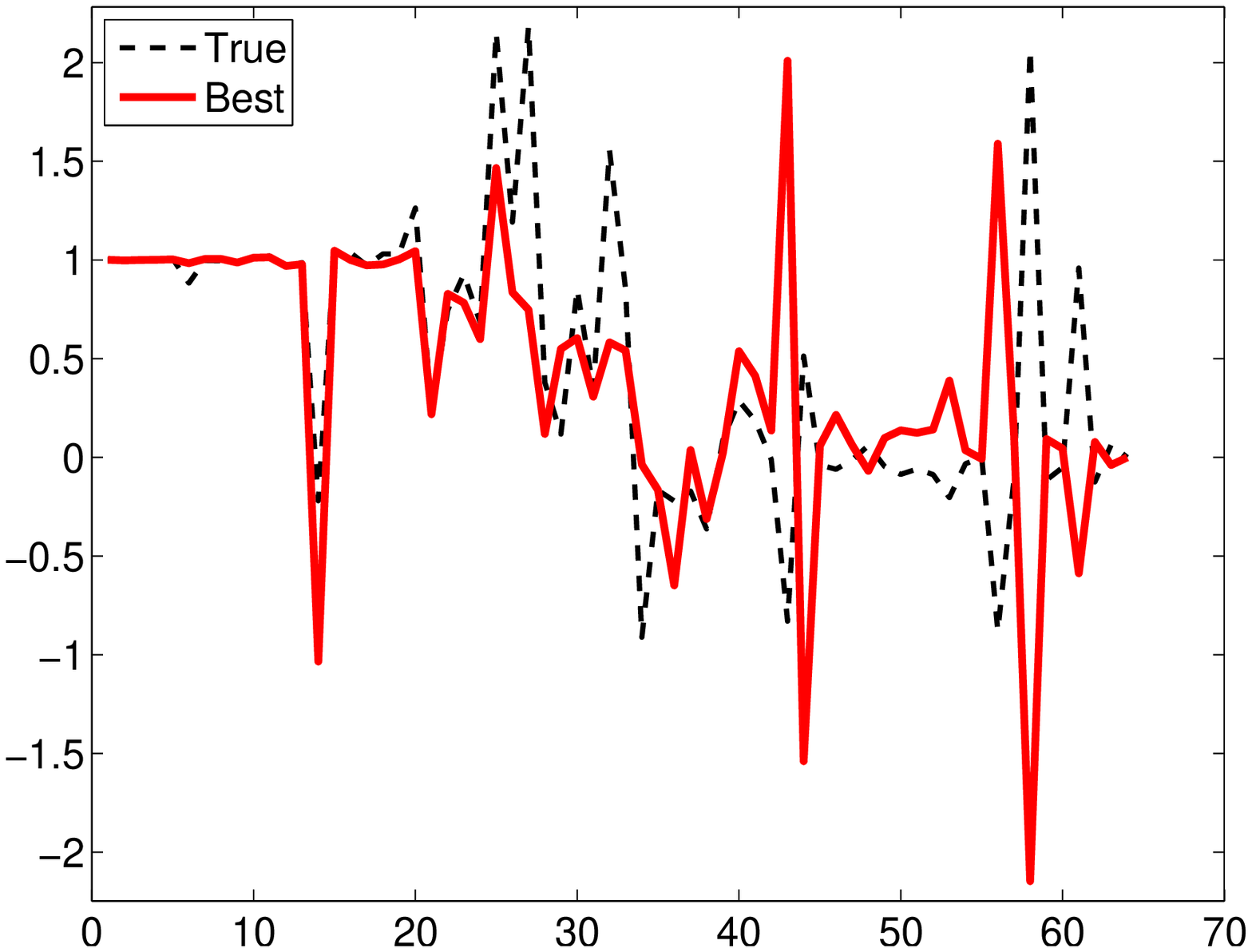}\\
SD\_P & ISRA\_P & HMZ\_P
\end{tabular}
\caption{Comparison of the $\mybf{x}_{\text{true}}$ filter factors (black dashed) for the \texttt{heat} test problem with 
the ones generated by projected gradient methods with the ISRA and HMZ scaling matrices at the iteration corresponding 
to the minimum error (red solid). The nonscaled filters are also reported (SD\_P).}
\label{fig:FiltriHeatProj}
\end{center}
\end{figure}

\noindent From the reconstruction errors and the filter factors we can make the following remarks:
\begin{itemize}
\item as expected, the presence of a scaling matrix provides an acceleration of the algorithms, resulting in a general lower 
number of iterations required to achieve the best reconstruction. In particular, the well-known fast convergence of the CGLS 
method is attested by the very few iterations needed;
\item the slower convergence seems to help the diagonal scalings ISRA and HMZ in providing better regularized solutions. 
The improvements in the performances are clearly visible also in the plots of the filter factors. 
The CGLS filters $\phi_i^{k+1}$ ($i=1,\ldots,n$) at each iteration preserve the ``regularity'' with respect to the index $i$ 
already observed in Figure \ref{fig:FiltriHeatNoScal} for the nonscaled methods. On the contrary, the ISRA and HMZ scalings 
succeed in reconstructing more faithfully the irregular filters profile of the true solution;
\item the effect of the scaling matrix, however, seems to provide some general improvements in the reconstruction errors 
also for the CGLS matrix, even if the differences with the nonscaled case are minimal. 
From the plots shown in Figure \ref{fig:FiltriHeatScal}, we can observe that these improvements are due to the 
slightly better reconstructions of the very first filters, corresponding to the higher singular values.
\end{itemize}
The same test problem has been used for the analysis of the projected algorithms, since its solution is a non-negative vector. 
In this case, we considered only the nonscaled method and the diagonal scaling matrices ISRA and HMZ, for which the projection 
is trivial (see Section \ref{secproj}). Moreover, we did not find any difference between the projection schemes \eqref{gradprojD} 
and \eqref{gradproj2D}, therefore we will present only the results obtained in the latter case. The best reconstruction errors 
and the corresponding numbers of iterations required by the different algorithms are reported again in Table \ref{Tab2}, while 
the filter factors are plotted in Figure \ref{fig:FiltriHeatProj}. In all cases, the presence of the constraint on the 
unknown produces a positive effect on the solution, attested by the reported lower reconstruction error. In particular, in the 
nonscaled case we can observe the irregular filter profile already remarked for the scaled methods, in agreement with the fact 
that also the presence of a matrix $\mybf{D}_k$ makes any filter factor related to all the singular values and singular vectors 
of the matrix $\mybf{A}$.

\begin{figure}[ht]
\begin{center}
 \begin{tabular}{ccc}
\includegraphics[width = 0.25\textwidth]{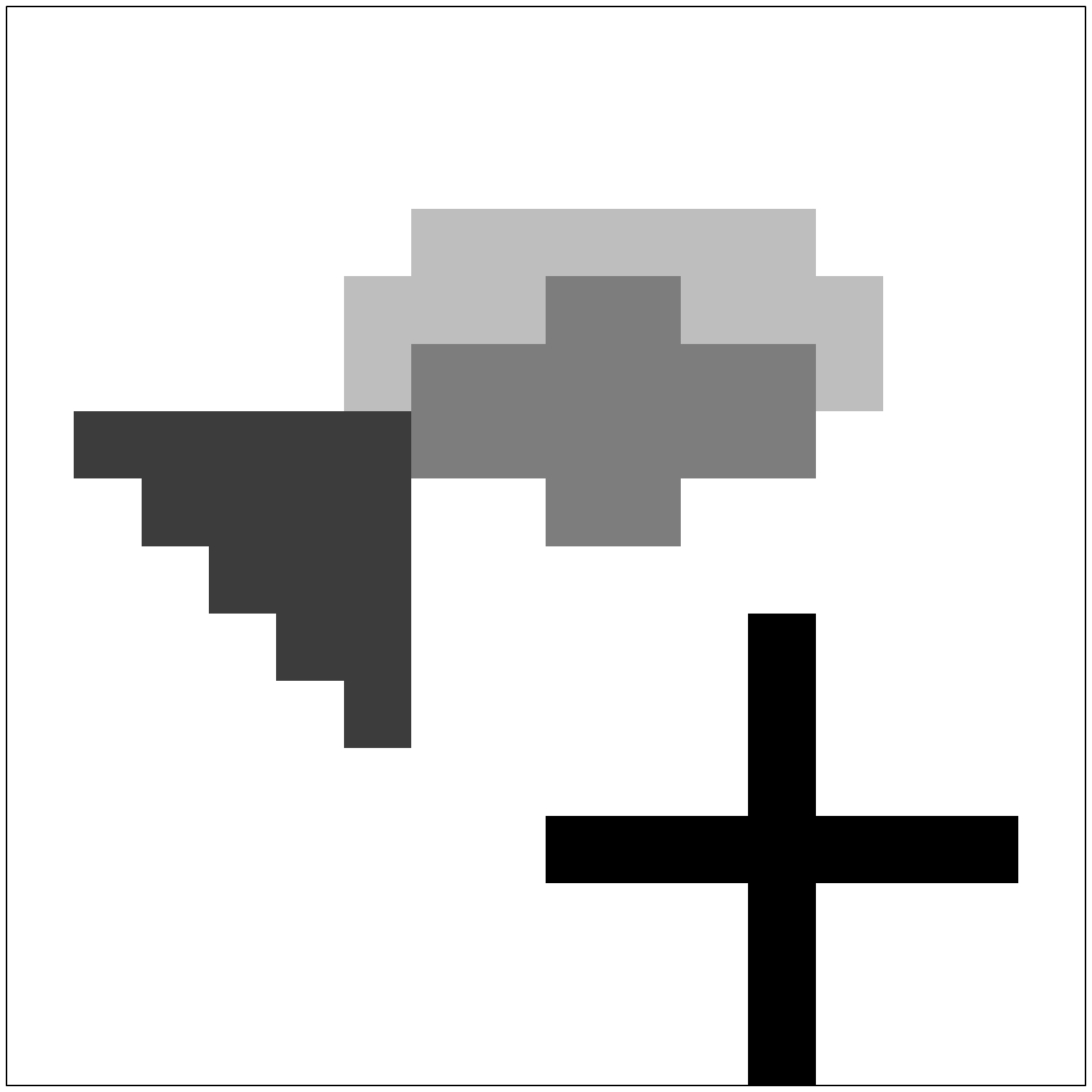}&
\includegraphics[width = 0.25\textwidth]{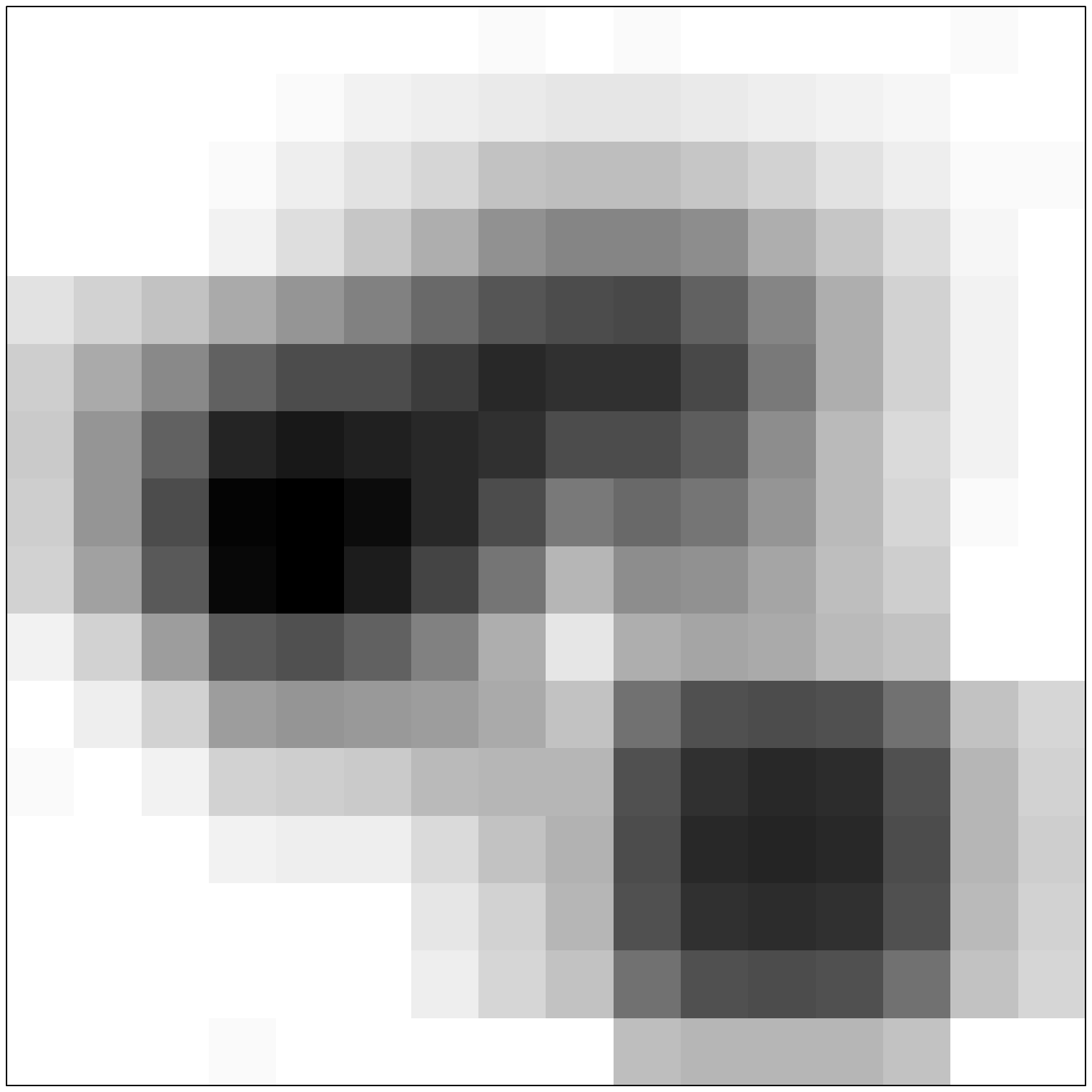}&
\includegraphics[width = 0.25\textwidth]{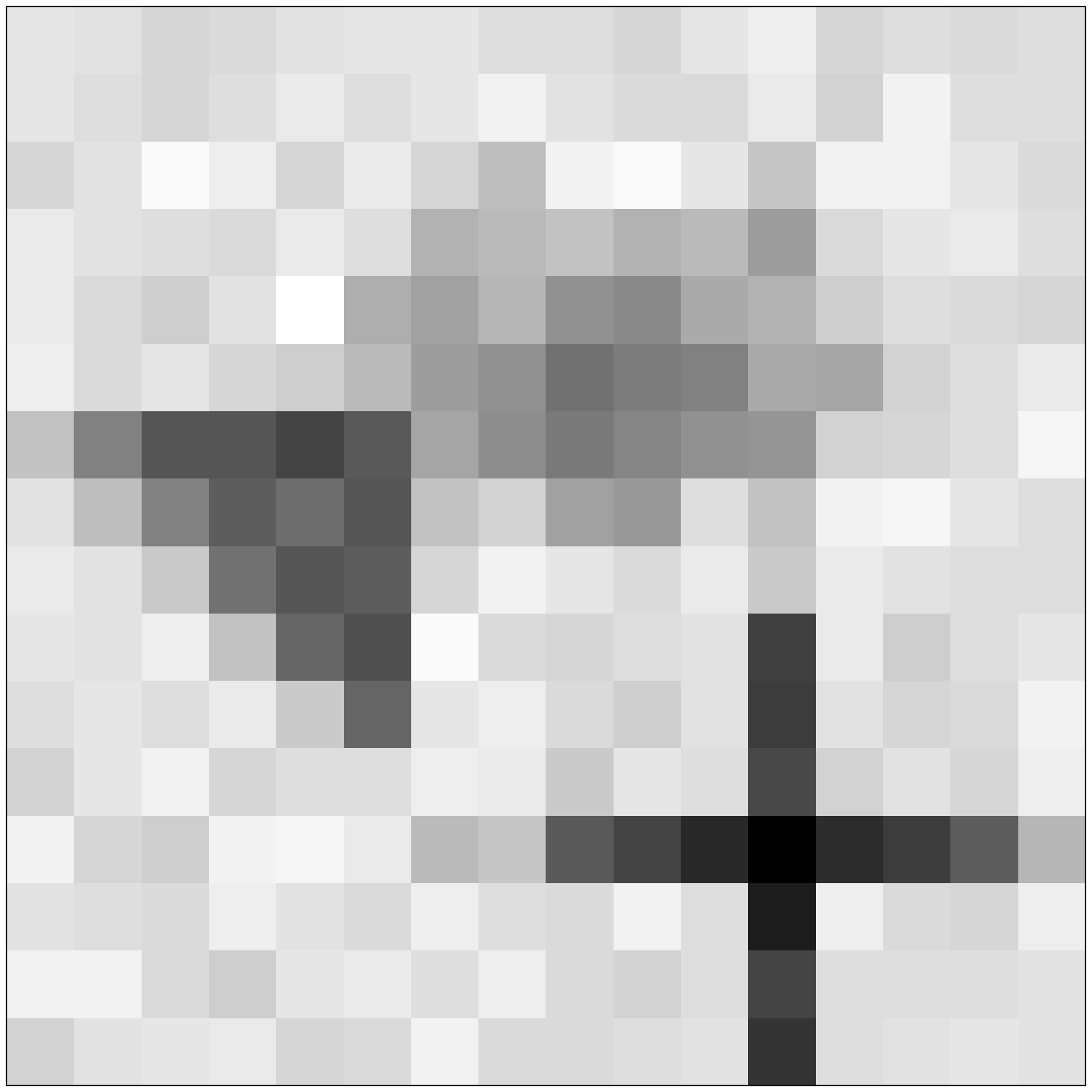}\\
True image & Measured image & SD \\
\includegraphics[width = 0.25\textwidth]{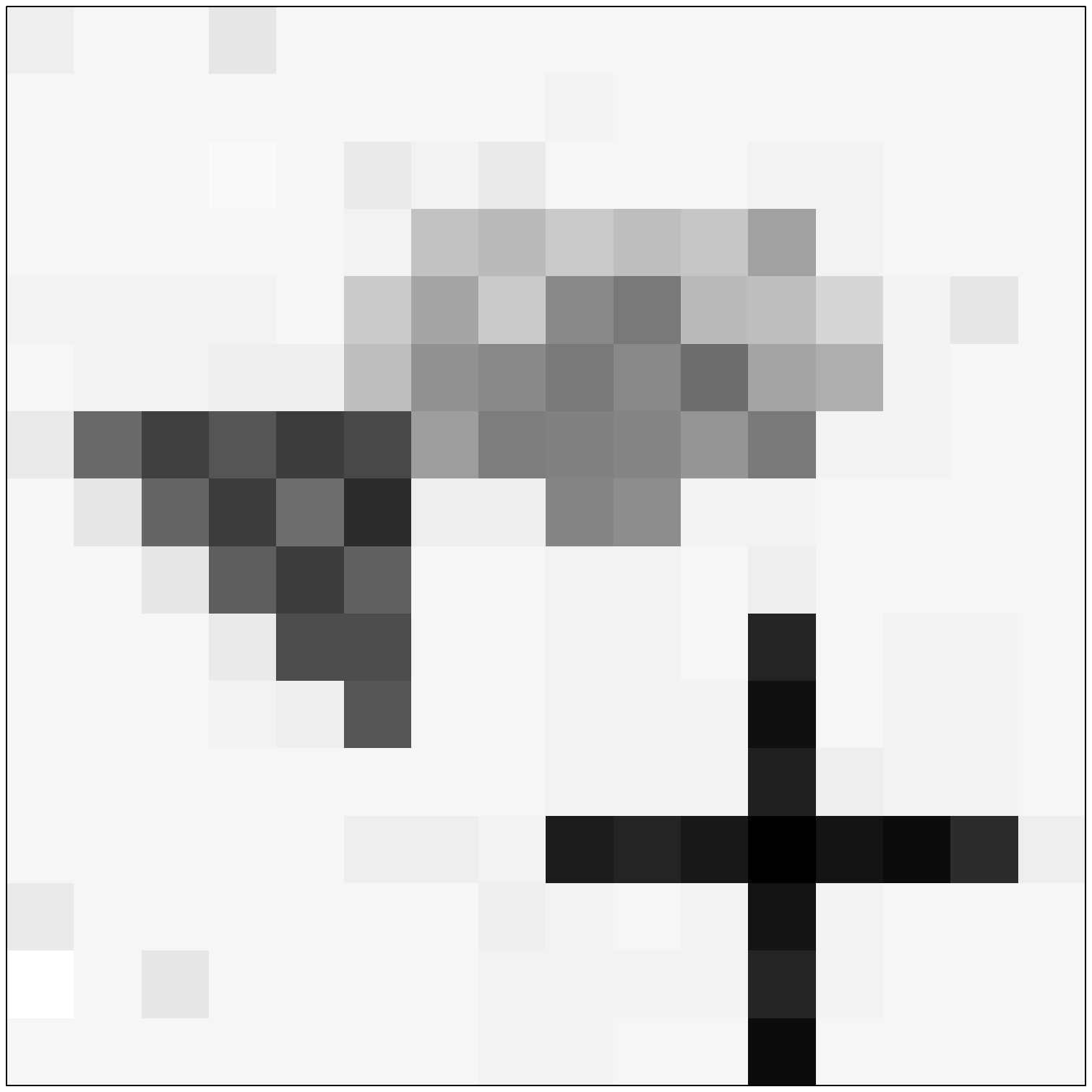}&
\includegraphics[width = 0.25\textwidth]{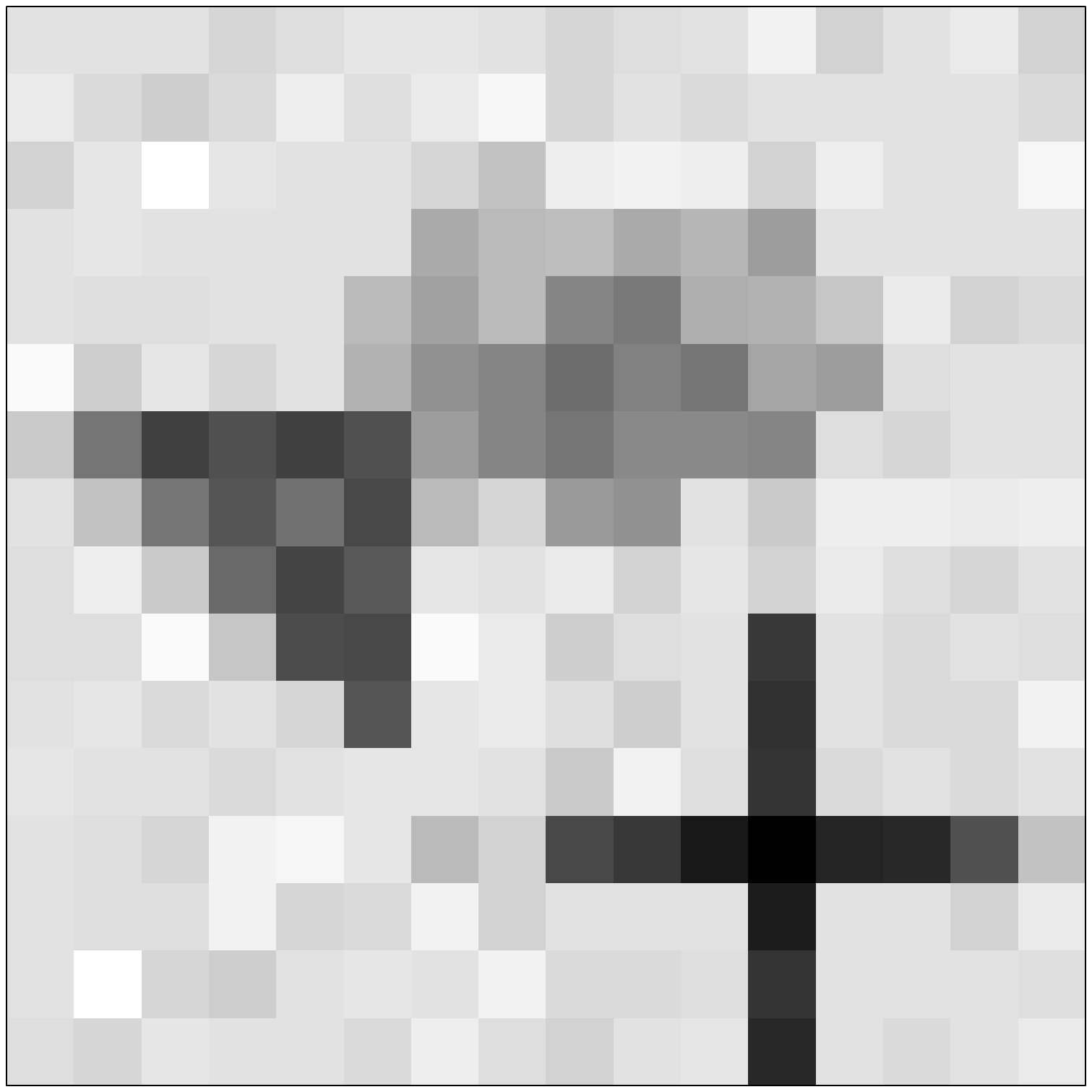}&
\includegraphics[width = 0.25\textwidth]{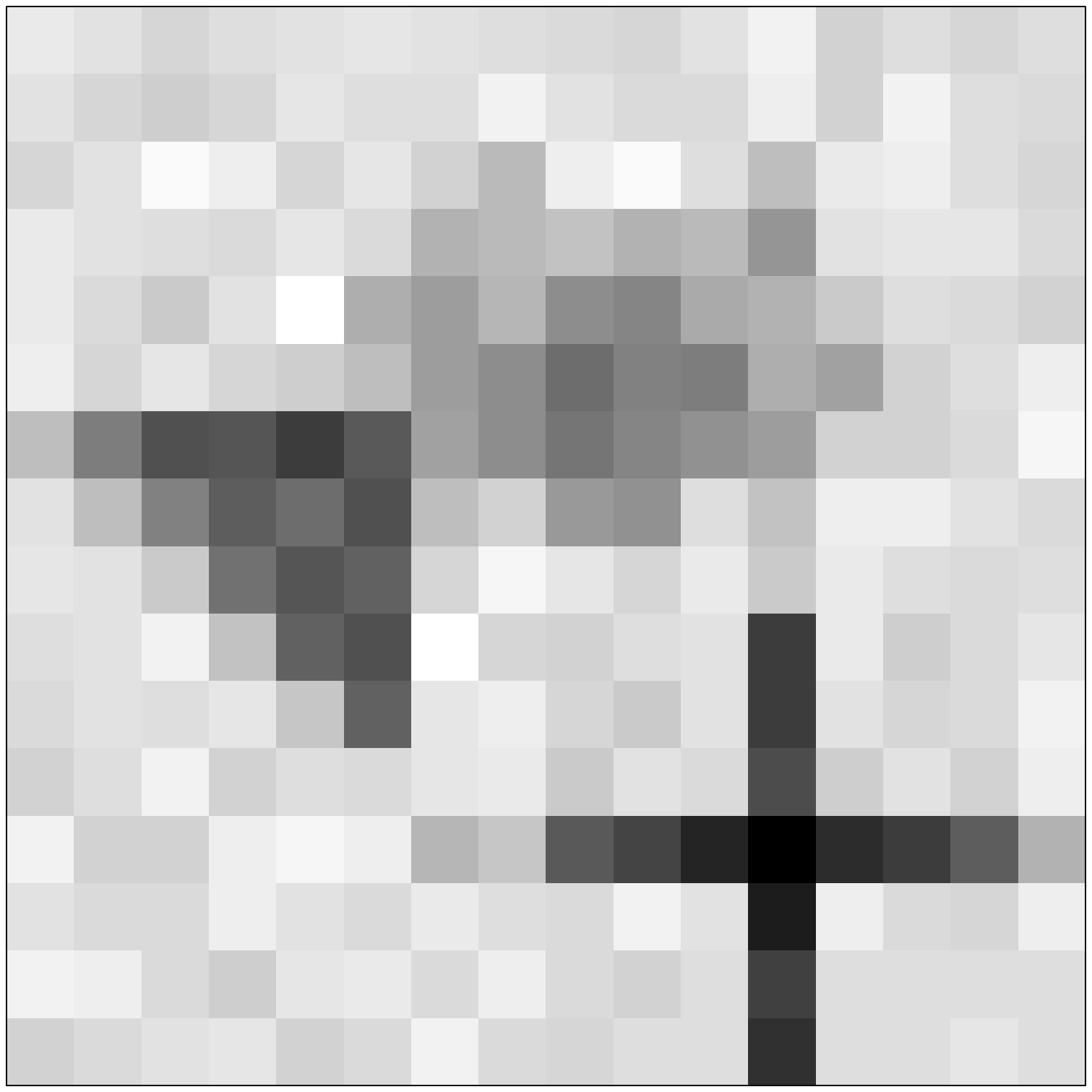}\\
ISRA & HMZ & CGLS \\
\includegraphics[width = 0.25\textwidth]{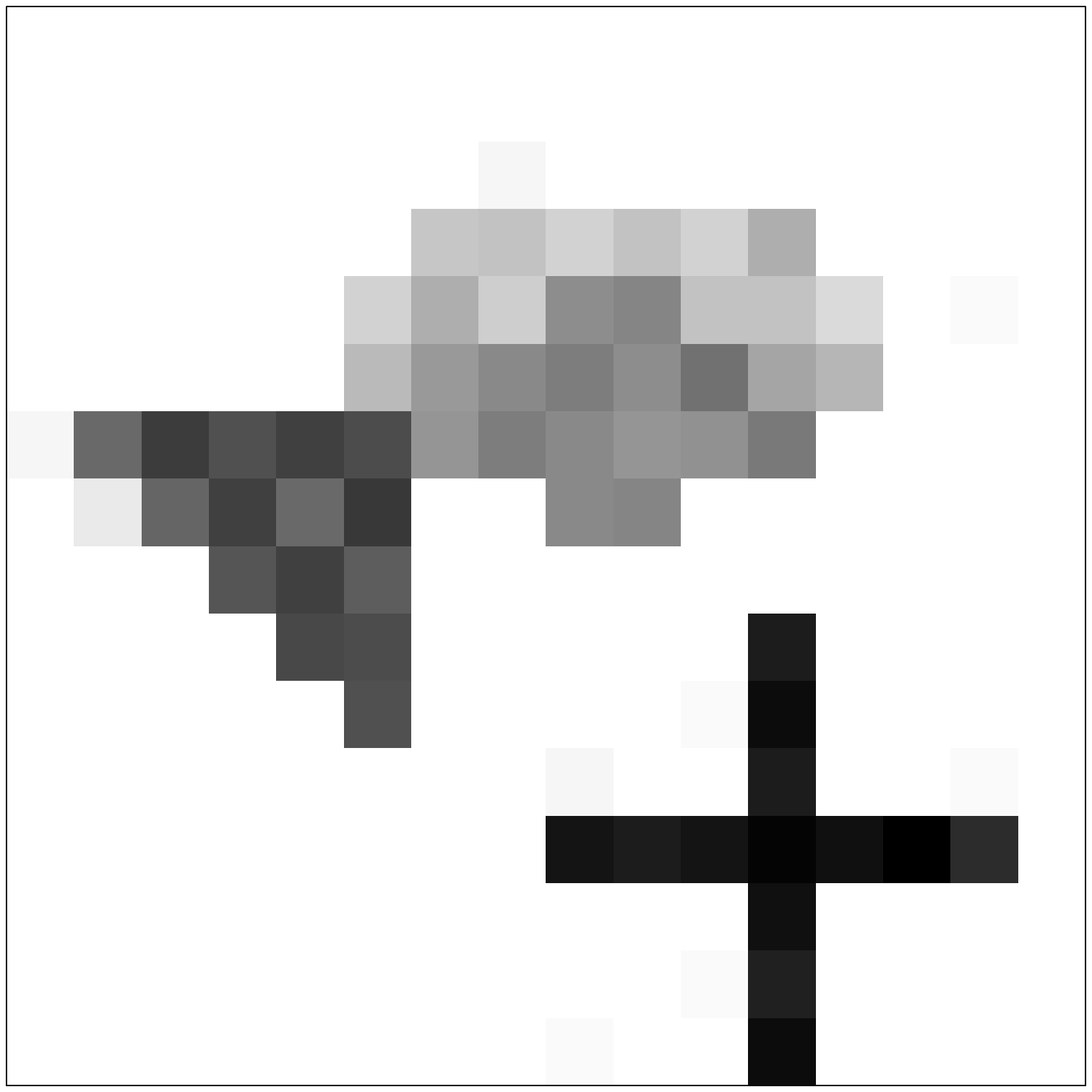}& 
\includegraphics[width = 0.25\textwidth]{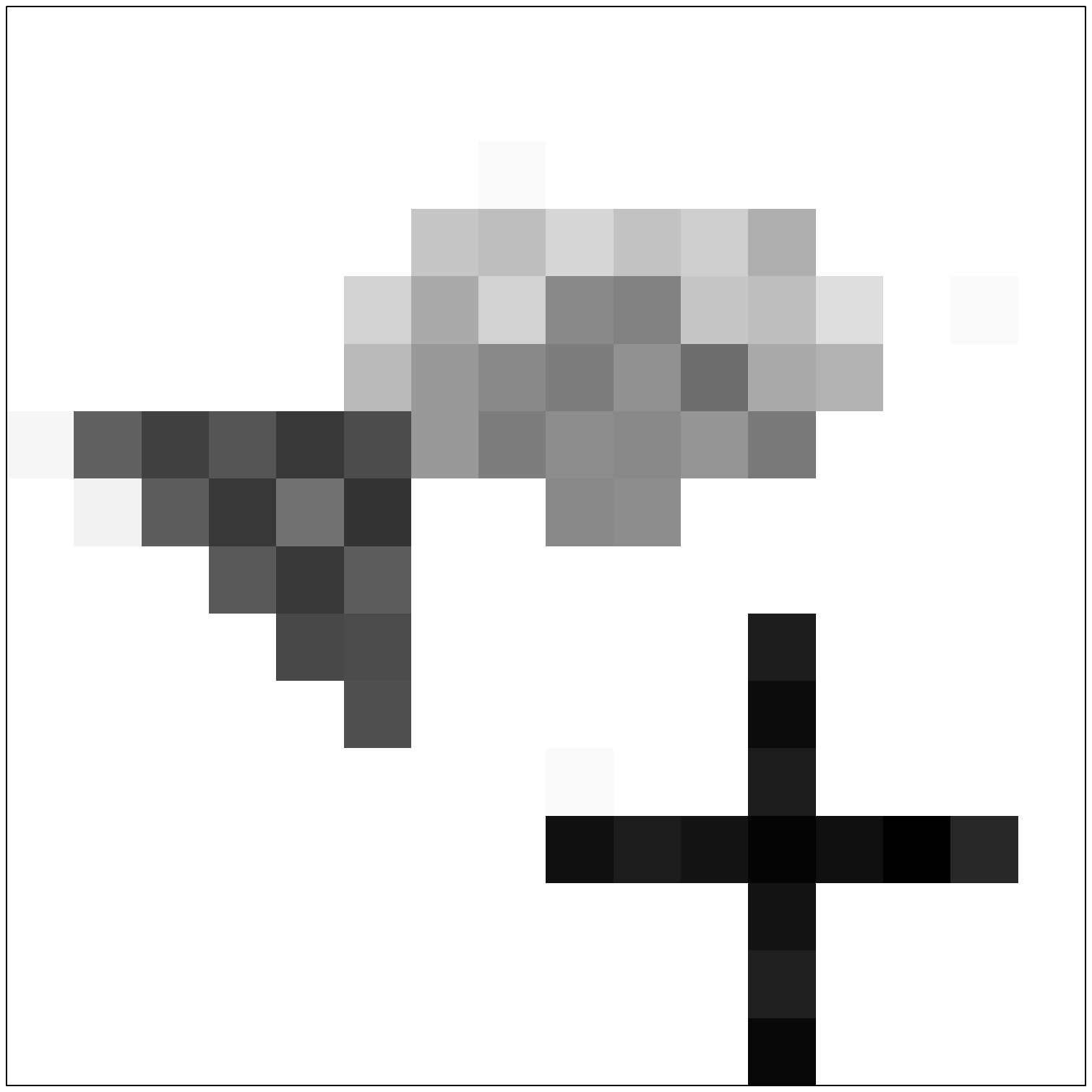}&
\includegraphics[width = 0.25\textwidth]{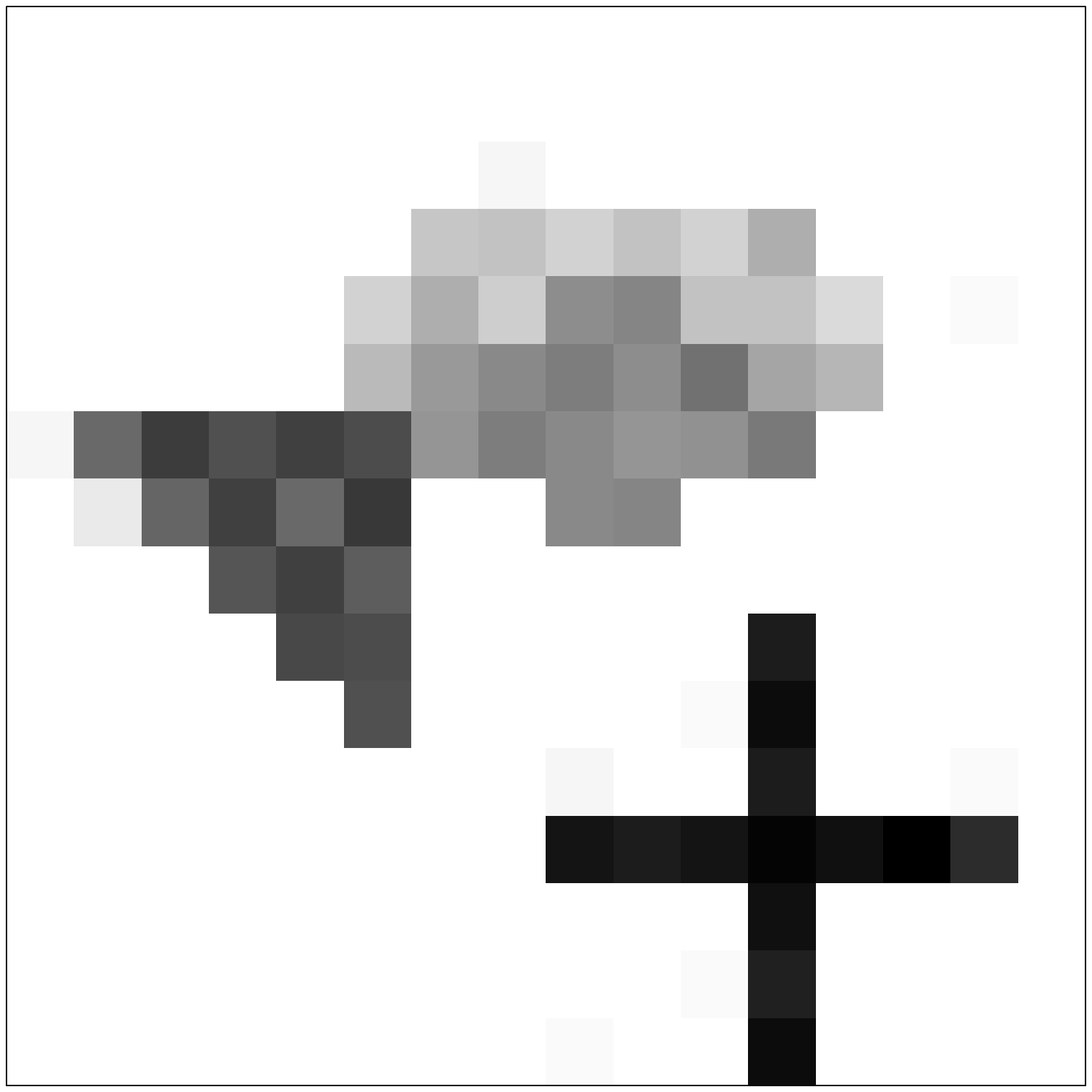}\\
SD\_P & ISRA\_P & HMZ\_P
 \end{tabular}
\caption{Comparison of the best reconstructed images for the \texttt{blur} test problem generated by gradient methods with 
different scaling matrices. The reconstruction with a nonscaled method is also reported (SD). The suffix `\_P' denotes the 
projected algorithms. True object and corresponding blurred noisy image are also showed.}\label{figblur1}
\end{center}
\end{figure}

\noindent Very similar considerations hold true also in the second test we carried out and concerning a 2-dimensional dataset. 
In particular, we used the \texttt{blur} test problem, taken again from Hansen's \textit{Regularization Tools} and simulating 
the degrading effect on a real image due to the action of a general acquisition system (see Figure \ref{figblur1}). 
Here we set $n = 256$, thus obtaining the true and measured images of $16 \times 16$ pixels. Moreover, we assumed the same 
parameters for the noise that we chose in the \texttt{heat} case.\\
A clear figure with the plots of the filter factors for all the considered gradient methods is hard to be produced, due to 
the irregularity of the filter profiles and the higher number of filters themselves with respect to the 1-dimensional test. 
However, we noticed again the regular behaviour in the unconstrained problem of the nondiagonal scalings, while the filter 
factors provided by ISRA and HMZ reproduced more precisely the ones computed with the true image. With the addition of the 
projection on the non-negative orthant, the reconstruction errors of both the nonscaled and the diagonally scaled algorithms 
decrease, reflecting the fact that a more faithful representation of the filters leads to a better reconstruction of the true image. 
In order to appreciate the positive effects on the regularized solutions, for this 2-dimensional case we reported in Figure 
\ref{figblur1} the restored images by the different algorithms. In particular, we can see that the gradient method with the 
ISRA scaling is able to remove most artifacts also in the unconstrained case. When the projected algorithms are used, instead, 
the presence of a scaling matrix seems to provide some benefits only in the reduced number of iterations required, 
since the best reconstruction error is already achievable with the nonscaled method.

\begin{table}
\begin{center}
\caption{Optimal numbers of iterations required by gradient methods with SD steplength and different scaling matrices 
for solving \texttt{blur}, with the corresponding minimum error. The suffix `\_P' denotes the projected algorithms.}
\label{Tabella2}
\begin{tabular}{l|cc}
\hline\noalign{\smallskip}
Method  & Iter & Min err \\
\noalign{\smallskip}\hline\noalign{\smallskip}
SD      & 1199 & 0.256 \\
CGLS    &   28 & 0.223 \\
ISRA    &  847 & 0.111 \\
HMZ     & 1352 & 0.165 \\
\noalign{\smallskip}\hline\noalign{\smallskip}
SD\_P   & 1870 & 0.089 \\
ISRA\_P & 1590 & 0.094 \\
HMZ\_P  &  707 & 0.089 \\
\noalign{\smallskip}\hline
\end{tabular}
\end{center}
\end{table}

\section{Conclusions}

In this paper we analyzed the regularizing effect of several gradient methods for the solution of the linear least-squares problem, 
i.e., the ability of the scheme to produce a sequence of vectors that, at a certain iteration, approximates as close as possible the 
true solution. The analysis we carried out has been made in terms of the ability of a given method to reproduce correctly the filter 
factors of the solution, accounting for the way in which the resulting sequence opposes the amplification effect of the noise on the 
data due to the presence of small singular values. The starting point of our work has been a paper of Nagy \& Palmer, in which the 
filter factors for nonscaled gradient methods have been formally calculated and numerically analyzed. In our paper we extended this 
analysis to the presence of scaling matrices in defining the descent direction, showing the analytical form of the corresponding filter 
factors and the advantages that can be obtained not only in terms of efficiency in decreasing the 
least-squares functional, but also in providing a better approximations of the unknown solution (as remarked e.g. in \cite{BON13}). 
Moreover, we also considered the case of non-negative constraints, and we generalized the expression of the filter factors to the 
projected gradient methods whose projection can be performed by means of the multiplication with a suitable diagonal matrix. 
We showed on some numerical examples the positive effect of the projection in reconstructing the irregular 
filter profiles of the true solution, in both scaled and nonscaled cases.

\section*{Acknowledgements}
This work has been partially supported by the Italian Spinner2013 PhD Project ``High-complexity inverse problems in biomedical applications and social systems'' and by a grant of the Italian Gruppo Nazionale per il Calcolo Scientifico (GNCS) - Istituto Nazionale di Alta Matematica (INdAM).

\end{document}